\documentclass[12pt]{article}
\usepackage{amsmath,amsfonts,amssymb,amsthm}
\usepackage{epsf,epsfig,longtable}

 \DeclareMathOperator{\Int}{Int}
\DeclareMathOperator{\Rep}{Rep} \DeclareMathOperator{\Per}{Per}

\DeclareMathOperator{\OrientSh}{OrientSh}
\DeclareMathOperator{\dist}{dist} \DeclareMathOperator{\diag}{diag}
\DeclareMathOperator{\StSh}{StSh} 
\DeclareMathOperator{\Cl}{Cl}
\DeclareMathOperator*{\KS}{KS}
\DeclareMathOperator{\D}{D}
\DeclareMathOperator{\OrbitSh}{OrbitSh}
 \DeclareMathOperator{\dd}{d}

\newcommand{\RR}{\ensuremath{\mathbb{R}}}

\newcommand{\Z}{\ensuremath{\mathbb{Z}}}
\newcommand{\CC}{\ensuremath{\mathbb{C}}}

\begin{document}

\newcommand{\al}{\mbox{$\alpha$}}
\newcommand{\si}{\mbox{$\sigma$}}
\newcommand{\de}{\mbox{$\delta$}}
\newcommand{\om}{\mbox{$\omega$}}
\newcommand{\De}{\mbox{$\Delta$}}
\newcommand{\ep}{\varepsilon}
\newcommand{\lam}{\mbox{$\lambda$}}
\newcommand{\La}{\mbox{$\Lambda$}}
\newcommand{\vp}{\mbox{$\varphi$}}
\newcommand{\gam}{\mbox{$\gamma$}}

\newcommand{\Ss}{\mbox{$\textbf{S}$}}
\newcommand{\ZZ}{\mbox{$\textbf{Z}$}}
\newcommand{\NN}{\mbox{$\textbf{N}$}}
\newcommand{\Cc}{\mbox{${\bf C}$}}
\newcommand{\Ff}{\mbox{${\cal F}$}}
\newcommand{\Bb}{\mbox{${\cal B}$}}
\newcommand{\Kk}{\mbox{${\cal K}$}}
\newcommand{\HH}{\mbox{${\cal H}$}}
\newcommand{\Tt}{\mbox{${\cal T}$}}

\newcommand{\KK}{\mbox{${\bf K}$}}

\newcommand{\Ws}{\mbox{$W^s$}}
\newcommand{\Wu}{\mbox{$W^u$}}
\newcommand{\Cone}{\mbox{$\Cc^1$}}
\newcommand{\ddt}{\mbox{$\frac{\dd }{\dd t}$}}

\newcommand{\sref}[1]{(\ref{#1})}
\newcommand{\setg}{\{g(t): \; t \in \RR\}}



\title{Vector Fields with the Oriented Shadowing Property}

\author{Sergei Yu.\ Pilyugin \footnote{Faculty of Mathematics and Mechanics, St.Petersburg State
University, University av.~28, 198504, St. Petersburg, Russia,
sp@sp1196.spb.edu}, Sergey B.\ Tikhomirov \footnote{Department of
Mathematics, National Taiwan University, No. 1, Section 4, Roosevelt
Road, Taipei 106, Taiwan, , sergey.tikhomirov@gmail.com.}\footnote{
Research of the second author is supported by NSC (Taiwan)
98-2811-M-002-061.}}

\maketitle

\begin{abstract} We give a description of the $\Cone$-interior
($\Int^1(\OrientSh)$) of the set of smooth vector fields on a smooth
closed manifold that have the oriented shadowing property. A special
class $\Bb$ of vector fields that are not structurally stable is
introduced. It is shown that the set $\Int^1(\OrientSh\setminus\Bb)$
coincides with the set of structurally stable vector fields. An
example of a field of the class $\Bb$ belonging to
$\Int^1(\OrientSh)$ is given. Bibliography: 18 titles.
\end{abstract}



\section{Introduction}

The theory of shadowing of approximate trajectories
(pseudotrajectories) in dynamical systems is now well developed
(see, for example, the monographs \cite{PilBook, Palm}). At the same
time, the problem of complete description of systems having the
shadowing property seems unsolvable. We have no hope to characterize
systems with the shadowing property in terms of the theory of
structural stability (such as hyperbolicity and transversality)
since the shadowing property is preserved under homeomorphisms of
the phase space (at least in the compact case), while the
above-mentioned properties are not.

The situation changes completely when we pass from the set of smooth
dynamical systems having the shadowing property (or some of its
analogs) to its $\Cone$-interior. It was shown by Sakai \cite{Sak}
that the $\Cone$-interior of the set of diffeomorphisms with the
shadowing property coincides with the set of structurally stable
diffeomorphisms. Later, a similar result was obtained for the set of
diffeomorphisms with the orbital shadowing property
\cite{PilRodSak}.

In this context, there is a real difference between the cases of
discrete dynamical systems generated by diffeomorphisms and systems
with continuous time (flows) generated by smooth vector fields. This
difference is due to the necessity of reparametrizing shadowing
trajectories in the latter case. One of the main goals of the
present paper is to show that this difference is crucial, and the
results for flows are essentially different from those for
diffeomorphisms.

Let us pass to the main definitions and results. Let $M$ be a smooth
closed (i.e., compact and boundaryless) manifold with Riemannian
metric $\dist$ and let $n = \dim M$. Consider a smooth ($\Cone$)
vector field on $X$ and denote by $\phi$ the flow of $X$. We denote
by
$$
O(x,\phi)=\{\phi(t,x):\;t\in\RR\}
$$
the trajectory of a point $x$ in the flow $\phi$; $O^+(x,\phi)$ and
$O^-(x,\phi)$ are the positive and negative semitrajectories,
respectively.

Fix a number $d>0$. We say that a mapping $g:\RR\to M$ (not
necessarily continuous) is a $d$-pseudotrajectory (both for the
field $X$ and flow $\phi$) if
\begin{equation}
\label{01} \dist(g(\tau+t),\phi(t,g(\tau)))<d\quad\mbox{for}\quad
\tau\in\RR,\;t\in[0,1].
\end{equation}

A reparametrization is an increasing homeomorphism $h$ of the line
$\RR$; we denote by $\Rep$ the set of all reparametrizations.

For $a>0$, we denote
$$
\Rep(a)=\left\{h\in\Rep:\;\left|\frac{h(t)-h(s)}{t-s}-1\right|<a,\quad
t,s\in\RR,\;t\neq s\right\}.
$$

In this paper, we consider the following three shadowing properties
(and the corresponding sets of dynamical systems).

We say that a vector field $X$ has the standard shadowing property
($X\in\StSh$) if for any $\ep>0$ we can find $d>0$ such that for any
$d$-pseudotrajectory $g(t)$ of $X$ there exists a point $p\in M$ and
a reparametrization $h\in\Rep(\ep)$ such that
\begin{equation}
\label{02} \dist(g(t),\phi(h(t),p))<\ep\quad\mbox{for}\quad t\in\RR.
\end{equation}

We say that a vector field $X$ has the oriented shadowing property
($X\in\OrientSh$) if for any $\ep>0$ we can find $d>0$ such that for
any $d$-pseudotrajectory of $X$ there exists a point $p\in M$ and a
reparametrization $h\in\Rep$ such that inequalities (\ref{02}) hold
(thus, it is not assumed that the reparametrization $h$ is close to
identity).

Finally, we say that a vector field $X$ has the orbital shadowing
property ($X\in\OrbitSh$) if for any $\ep>0$ we can find $d>0$ such
that for any $d$-pseudotrajectory of $X$ there exists a point $p\in
M$ such that
$$
\dist_H(\Cl O(p,\phi),\Cl \{g(t):\;t\in\RR\})<\ep,
$$
where $\dist_H$ is the Hausdorff distance.

Let us note that the standard shadowing property is equivalent to
the strong pseudo orbit tracing property (POTP) in the sense of
Komuro \cite{Kom}; the oriented shadowing property was called the
normal POTP by Komuro \cite{Kom} and the POTP for flows by Thomas
\cite{Tho}.

We consider the following $\Cone$ metric on the space of smooth
vector fields: If $X$ and $Y$ are vector fields of class $\Cone$, we
set
$$
\rho_1(X,Y)=\max_{x\in M}\left(|X(x)-Y(x)|+ \left\|\frac{\partial
X}{\partial x}(x)-\frac{\partial Y}{\partial x}(x) \right\|\right),
$$
where $|.|$ is the norm on the tangent space $T_xM$ generated by the
Riemannian metric $\mbox{dist}$, and $\|.\|$ is the corresponding
operator norm for matrices.

For a set $A$ of vector fields, $\Int^1(A)$ denotes the interior of
$A$ in the $\Cone$ topology generated by the metric $\rho_1$.

Let us denote by $\Ss$ and $\NN$ the sets of structurally stable and
nonsingular vector fields, respectively.

The only result in the problem under study was recently published by
Lee and Sakai \cite{LeeSak}: $\Int^1(\StSh\cap\;\NN)\subset\Ss$.

To formulate our main results, we need one more definition.

Let us say that a vector field $X$ belongs to the class $\Bb$ if $X$
has two hyperbolic rest points $p$ and $q$ (not necessarily
different) with the following properties:

(1) The Jacobi matrix $DX(q)$ has two complex conjugate eigenvalues
$\mu_{1,2}=a_1\pm ib_1$ of multiplicity one with $a_1<0$ such that
if $\lam\neq\mu_{1,2}$ is an eigenvalue of $DX(q)$ with
$\mbox{Re}\lam<0$, then $\mbox{Re}\lam<a_1$;

(2) the Jacobi matrix $DX(p)$ has two complex conjugate eigenvalues
$\nu_{1,2}=a_2\pm ib_2$ with $a_2>0$ of multiplicity one such that
if $\lam\neq\nu_{1,2}$ is an eigenvalue of $DX(p)$ with
$\mbox{Re}\lam>0$, then $\mbox{Re}\lam>a_2$;

(3) the stable manifold $\Ws(p)$ and the unstable manifold $\Wu(q)$
have a trajectory of nontransverse intersection.

Condition (1) above means that the ``weakest" contraction in
$\Ws(q)$ is due to the eigenvalues $\mu_{1,2}$ (condition (2) has a
similar meaning).
\medskip

{\bf Theorem 1. } $\Int^1(\OrientSh\setminus\Bb)=\Ss.$
\medskip

Let us note that Theorem 1 was stated (without a proof) in the
author's short note \cite{PilTikh}. Let us also note that if
$\mbox{dim}M\leq 3$, then $\Int^1(\OrientSh)=\Ss$ (which also was
stated in \cite{PilTikh} and proved by the second author in
\cite{Tikh}; in \cite{Tikh}, it was also shown that if LipSh is the
set of vector fields that have an analog of the standard shadowing
property with $\ep$ replaced by $Ld$, then
$\Int^1(\mbox{LipSh})=\Ss$).
\medskip

{\bf Theorem 2. } $\Int^1(\OrientSh)\cap\Bb\neq\emptyset.$
\medskip

{\bf Theorem 3. } $\Int^1(\OrbitSh\cap\;\NN)\subset\Ss.$
\medskip

Let us note that Theorem 3 generalizes the above-mentioned result by
Lee and Sakai.

The structure of the paper is as follows: In Sec. 2, we prove
Theorem 1 and discuss the proof of Theorem 3; in Sec. 3, we prove
Theorem 2.

\section{Proof of Theorem 1}

First we introduce some notation.

We denote by $B(a,A)$ the $a$-neighborhood of a set $A\subset M$.

The term ``transverse section" will mean a smooth open disk in $M$
of codimension 1 that is transverse to the flow $\phi$ at any of its
points.

Let $\Per(X)$ denote the set of rest points and closed orbits of a
vector field $X$.

Let us recall that $X$ is called a Kupka-Smale field ($X\in\KS$) if

(KS1) any trajectory in $\Per(X)$ is hyperbolic;

(KS2) stable and unstable manifolds of trajectories from $\Per(X)$
are transverse.

The proof of Theorem 1 is based on the following result (see
\cite{Gan}): $\Int^1(\KS)=\Ss$.

Let ${\cal T}$ denote the set of vector fields $X$ that have
property (KS1). Our first lemma is applied in the proofs of both
Theorems 1 and 3; for this purpose, we formulate and prove it for
the set $\OrbitSh$.
\medskip

{\bf Lemma 1.}
\begin{equation}
 \label{l2}
\mbox{Int}^1(\mbox{OrbitSh})\subset{\cal T}.
\end{equation}

{\em Proof.} To get a contradiction, let us assume that that there
exists a vector field $X\in \mbox{Int}^1(\mbox{OrbitSh})$ that does
not have property (KS1), i.e., the set $\mbox{Per}(X)$ contains a
trajectory $p$ that is not hyperbolic.

Let us first consider the case where $p$ is a rest point. Identify
$M$ with $\RR^n$ in a neighborhood of $p$. Applying an arbitrarily
$\Cone$-small perturbation of the field $X$, we can find a field
$Y\in \mbox{Int}^1(\mbox{OrbitSh})$ that is linear in a neighborhood
$U$ of $p$ (we also assume that $p$ is the origin of $U$).

(Here and below in the proof of Lemma 1, all the perturbations are
$\Cone$-small perturbations that leave the field in
$\mbox{Int}^1(\mbox{OrbitSh})$; we denote the perturbed fields by
the same symbol $X$ and their flows by $\phi$.)

Then trajectories of $X$ in $U$ are governed by a differential
equation
\begin{equation}
 \label{l3}
\dot{x}=Px,
\end{equation}
where the matrix $P$ has an eigenvalue $\lambda$ with
$\mbox{Re}\lambda=0$.

Consider first the case where $\lambda=0$. We perturb the field $X$
(and change coordinates, if necessary) so that, in Eq. (\ref{l3}),
the matrix $P$ is block-diagonal,
\begin{equation}
 \label{l4}
P=\mbox{diag}(0,P_1),
\end{equation}
and $P_1$ is an $(n-1)\times(n-1)$ matrix.

Represent coordinate $x$ in $U$ as $x=(y,z)$ with respect to
(\ref{l4}); then
$$
\phi(t,(y,z))=\left(y,\exp(P_1t)z\right)
$$
in $U$.

Take $\epsilon>0$ such that $B(4\epsilon,p)\subset U$. To get a
contradiction, assume that $X\in\mbox{OrbitSh}$; let $d$ correspond
to the chosen $\epsilon$.

Fix a natural number $m$ and consider the following mapping from
$\RR$ into $U$:
$$
g(t)=\begin{cases}
y=-2\epsilon, \quad z=0;\quad t\leq 0,\\
y=-2\epsilon+t/m, \quad z=0;\quad 0<t<4m\epsilon,\\
y=2\epsilon, \quad z=0;\quad 4m\epsilon<t.
\end{cases}
$$

Since the mapping $g$ is continuous, piecewise differentiable, and
either $\dot{y}=0$ or $\dot{y}=1/m$, $g$ is a $d$-pseudotrajectory
for large $m$.

Any trajectory of $X$ in $U$ belongs to a plane $y=\mbox{const}$;
hence,
$$
\mbox{dist}_H(\mbox{Cl}(O(q,\phi)),\mbox{Cl}(\{g(t):t\in \RR\}))\geq
2\epsilon
$$
for any $q$. This completes the proof in the case considered.

Similar reasoning works if $p$ is a rest point and the matrix $P$ in
(\ref{l3}) has a pair of eigenvalues $\pm ib, b\neq 0$.

Now we assume that $p$ is a nonhyperbolic closed trajectory. In this
case, we perturb the vector field $X$ in a neighborhood of the
trajectory $p$ using the perturbation technique developed by Pugh
and Robinson in \cite{PughRob}. Let us formulate their result (which
will be used below several times).
\medskip

{\bf Pugh-Robinson pertubation. }{\em Assume that $r_1$ is not a
rest point of a vector field $X$. Let $r_2=\phi(\tau,r_1)$, where
$\tau>0$. Let $\Sigma_1$ and $\Sigma_2$ be two small transverse
sections such that $r_i\in\Sigma_i, i=1,2$. Let $\si$ be the local
Poincar\'e transformation generated by these transverse sections.

Consider a point $r'=\phi(\tau',r_1)$, where $\tau'\in(0,\tau)$, and
let $U$ be an arbitrary open set containing $r'$.

Fix an arbitrary $C^1$-neighborhood $F$ of the field $X$.

There exist positive numbers $\ep_0$ and $\De_0$ with the following
property: if $\si'$ is a local diffeomorphism from the
$\De_0$-neighborhood of $r_1$ in $\Sigma_1$ into $\Sigma_2$ such
that}
$$
\mbox{dist}_{C^1}(\si,\si')<\ep_0,
$$
{\em then there exists a vector field $X'\in F$ such that

(1) $X'=X$ outside $U$;

(2) $\si'$ is the local  Poincar\'e transformation generated by the
sections $\Sigma_1$ and $\Sigma_2$ and trajectories of the field}
$X'$.
\medskip

Let $\om$ be the least positive period of the nonhyperbolic closed
trajectory $p$. We fix a point $\pi\in p$, local coordinates in
which $\pi$ is the center, and a hyperplane $\Sigma$ of codimension
1 transverse to the vector $F(\pi)$. Let $y$ be coordinate in
$\Sigma$.

Let $\si$ be the local Poincar\'e transformation generated by the
transverse section $\Sigma$; denote $P=D\si(0)$. Our assumption
implies that the matrix $P$ is not hyperbolic. In an arbitrarily
small neighborhood of the matrix $P$, we can find a matrix $P'$ such
that $P'$ either has a real eigenvalue with unit absolute value of
multiplicity 1 or a pair of complex conjugate eigenvalues with unit
absolute value of multiplicity 1. In both cases, we can choose
coordinates $y=(v,w)$ in $\Sigma$ in which
\begin{equation}
 \label{l5}
P'=\mbox{diag}(Q,P_1),
\end{equation}
where $Q$ is a $1\times 1$ or $2\times 2$ matrix such that
$|Qv|=|v|$ for any $v$.

Now we can apply the Pugh-Robinson perturbation (taking
$r_1=r_2=\pi$ and $\Sigma_1=\Sigma_2=\Sigma$) that modifies $X$ in a
small neighborhood of the point $\phi(\om/2,\pi)$ and such that, for
the perturbed vector field $X'$, the local Poincar\'e transformation
generated by the transverse section $\Sigma$ is given by $y\mapsto
P'y$.

Clearly, in this case, the trajectory of $\pi$ in the field $X'$ is
still closed (with some period $\om'$). As was mentioned, we assume
that $X'$ has the orbital shadowing property (and write $X,\phi,\om$
instead of $X',\phi',\om'$).

We introduce in a neighborhood of the point $\pi$ coordinates
$x=(x',y)$, where $x'$ is one-dimensional (with axis parallel to
$X(\pi)$), and $y$ has the above-mentioned property.

Of course, the new coordinates generate a new metric, but this new
metric is equivalent to the original one; thus, the corresponding
shadowing property (or its absence) is preserved.

We need below one more technical statement.
\medskip

{\bf LE (local estimate)}.{\em There exists a neighborhood $W$ of
the origin in $\Sigma$ and constants $l,\delta_0>0$ with the
following property: if $z_1\in \Sigma\cap W$ and
$|z_2-z_1|<\delta<\delta_0$, then we can represent $z_2$ as
$\phi(\tau,z'_2)$ with $z'_2\in\Sigma$ and}
\begin{equation}
 \label{l01}
|\tau|,\;|z'_2-z_1|<l\delta.
\end{equation}

This statement is an immediate corollary of the theorem on local
rectification of trajectories (see, for example, \cite{Arnold}): In
a neighborhood of a point that is not a rest point, the flow of a
vector field of class $C^1$ is diffeomorphic to the family of
parallel lines along which points move with unit speed (and it is
enough to note that a diffemorphic image of $\Sigma$ is a smooth
submanifold transverse to lines of the family).

We may assume that the neighborhood $W$ in LE is so small that for
$y\in \Sigma\cap W$, the function $\al(y)$ (the time of first return
to $\Sigma$) is defined, and that the point $\phi(\al(v,w),(0,v,w))$
has coordinates $(Qv,P_1w)$ in $\Sigma$.

Let us take a neighborhood $U$ of the trajectory $p$ such that if
$r\in U$, then the first point of intersection of the positive
semitrajectory of $r$ with $\Sigma$ belongs to $W$.

Take $a>0$ such that the $4a$-neighborhood of the origin in $\Sigma$
is a subset of $W$. Fix
$$
\epsilon<\min\left(\delta_0,\frac{a}{4l}\right),
$$
where $\delta_0$ and $l$ satisfy the LE. Let $d$ correspond to this
$\epsilon$ (in the definition of the orbital shadowing property).

Take $y_0=(v_0,0)$ with $|v_0|=a$. Fix a natural number $N$ and set
$$
\al_k=\al\left(\left(\frac{k}{N}Q^kv_0,0\right)\right),\quad
k\in[0,N-1),
$$
$$
\beta_0=0,\quad \beta_k=\al_1+\dots+\al_k,
$$
and
$$
g(t)=\begin{cases}
\phi(t,(0,0,0)),\quad t<0;\\
\phi\left(t-\beta_k,\left(0,\frac{k}{N}Q^kv_0,0\right)\right),\quad
\beta_k\leq
t<\beta_{k+1},\;k\in[0,N-1);\\
\phi\left(t-\beta_N,\left(0,Q^Nv_0,0\right)\right),\quad
t\geq \beta_N.\\
\end{cases}
$$

Note that for any point $y=(v,0)$ of intersection of the set
$\{g(t):t\in \RR\}$ with $\Sigma$, the inequality $|v|\leq a$ holds.
Hence, we can take $a$ so small that
$$
B(2a,\mbox{Cl}(\{g(t):t\in \RR\}))\subset U.
$$

Since
$$
\left|\frac{k}{N}Q^{k+1}v_0-\frac{k+1}{N}Q^{k+1}v_0\right|=
\frac{a}{N}\to 0, \quad N\to \infty,
$$
$g(t)$ is a $d$-pseudotrajectory for large $N$.

Assume that there exists a point $q$ such that
$$
\mbox{dist}_H(\mbox{Cl}(O(q,\phi)),\mbox{Cl}(\{g(t):t\in
\RR\}))<\epsilon.
$$
In this case, $O(q,\phi)\subset U$, and there exist points
$q_1,q_2\in O(q,\phi)$ such that
$$
|q_1|=|q_1-(0,0,0)|<\epsilon
$$
and
$$
|q_2-(0,Q^Nv_0,0)|<\epsilon.
$$
By the choice of $\epsilon$, there exist points $q'_1,q'_2\in
O(q,\phi)\cap\Sigma$ such that
$$
|q'_1|<l\epsilon<a/4\quad\mbox{and}\quad
|q'_2-Q^Nv_0|<l\epsilon<a/4.
$$
Let $q'_1=(0,v_1,w_1)$ and $q'_2=(0,v_2,w_2)$. Since these points
belong to the same trajectory that is contained in $U$,
$|v_1|=|v_2|$. At the same time,
$$
|v_1|<a/4, \quad |v_2-Q^Nv_0|<a/4, \quad\mbox{and}\quad |Q^Nv_0|=a,
$$
and we get a contradiction which proves our lemma.

To complete the proof of Theorem 1, we show that any vector field
$$
X\in\mbox{Int}^1(\mbox{OrientSh}\setminus{\cal B})
$$
has property (KS2).

To get a contradiction, let us assume that there exist trajectories
$p,q\in \mbox{Per}(X)$ for which the unstable manifold $W^u(q)$ and
the stable manifold $W^s(p)$ have a point $r$ of nontransverse
intersection. We have to consider separately the following two
cases.

Case (B1): $p$ and $q$ are rest points of the flow $\phi$.

Case (B2): either $p$ or $q$ is a closed trajectory.

Case (B1). Since $X\notin{\cal B}$, we may assume (after an
additional perturbation, if necessary) that the eigenvalues
$\lambda_1,\dots,\lambda_u$ with $\mbox{Re}\lambda_j>0$ of the
Jacobi matrix $DX(p)$ have the following property:
$$
\mbox{Re}\lambda_j>\lambda_1>0,\quad j=2,\dots,u
$$
(where $u$ is the dimension of $W^u(p)$). This property means that
there exists a one-dimensional ``direction of weakest expansion'' in
$W^u(p)$.

If this is not the case, then our assumption that $X\notin{\cal B}$
implies that the eigenvalues $\mu_1,\dots,\mu_s$ with
$\mbox{Re}\mu_j<0$ of the Jacobi matrix $DX(q)$ have the following
property:
$$
\mbox{Re}\mu_j<\mu_1<0,\quad j=2,\dots,s
$$
(where $s$ is the dimension of $W^s(q)$). If this condition holds,
we reduce the problem to the previous case by passing from the field
$X$ to the field $-X$ (clearly, the fields $X$ and $-X$ have the
oriented shadowing property simultaneously).

Making a perturbation (in this part of the proof, we always assume
that the perturbed field belongs to the set
$\mbox{OrientSh}\setminus{\cal B}$), we may ``linearize'' the field
$X$ in a neighborhood $U$ of the point $p$; thus, trajectories of
$X$ in $U$ are governed by a differential equation
$$
\dot{x}=Px,
$$
where
\begin{equation}
 \label{l9}
P=\mbox{diag}(P_s,P_u),\quad P_u=\mbox{diag}(\lambda,P_1),\quad
\lambda>0,
\end{equation}
$P_1$ is a $(u-1)\times(u-1)$ matrix for which there exist constants
$K>0$ and $\mu>\lambda$ such that
\begin{equation}
 \label{l09}
\|\exp(-P_1t)\|\leq K^{-1}\exp(-\mu t),\quad t\geq 0,
\end{equation}
and $\mbox{Re}\lambda_j<0$ for the eigenvalues $\lambda_j$ of the
matrix $P_s$.

Let us explain how to perform the above-mentioned perturbations
preserving the nontransversalty of $W^u(q)$ and $W^s(p)$ at the
point $r$ (we note that a similar reasoning can be used in
``replacement" of a component of intersection of $W^u(q)$ with a
transverse section $\Sigma$ by an affine space, see the text
preceding Lemma 2 below).

Consider points $r^*=\phi(\tau,r)$, where $\tau>0$, and
$r'=\phi(\tau',r)$, where $\tau'\in(0,\tau)$. Let $\Sigma$ and
$\Sigma^*$ be small transverse sections that contain the points $r$
and $r^*$. Take small neighborhoods $V$ and $U'$ of $p$ and $r'$,
respectively, so that the set $V$ does not intersect the ``tube''
formed by pieces of trajectories through points of $U'$ whose
endpoints belong to $\Sigma$ and $\Sigma^*$. In this case, if we
perturb the vector field $X$ in $V$ and apply the Pugh-Robinson
perturbation in $U'$, these perturbations are ``independent.''

We perturb the vector field $X$ in $V$ obtaining vector fields $X'$
that are linear in small neighborhoods $V'\subset V$ and such that
the values $\rho_1(X,X')$ are arbitrarily small.

Let $\gam_s$ and $\gam^*_s$ be the components of intersection of the
stable manifold $W^s(p)$ (for the field $X$) with $\Sigma$ and
$\Sigma^*$ that contain the points $r$ and $r^*$, respectively.

Since the stable manifold of a hyperbolic rest point depends (on its
compact subsets) $C^1$-smoothly on $C^1$-small perturbations, the
stable manifolds $W^s(p)$ (for the perturbed fields $X'$) contain
components $\gam'_s$ of intersection with $\Sigma^*$ that converge
(in the $C^1$ metric) to $\gam^*_s$.

Now we apply the Pugh-Robinson perturbation in $U'$ and find a field
$X'$ in an arbitrary $C^1$ neighborhood of $X$ such that the local
Poincar\'e transformation generated by the field $X'$ and sections
$\Sigma$ and $\Sigma^*$ takes $\gam'_s$ to $\gam_s$ (which means
that the nontransversality at $r$ is preserved).

We introduce in $U$ coordinates $x=(y;v,w)$ according to (\ref{l9}):
$y$ is coordinate in the $s$-dimensional ``stable'' subspace
(denoted $E^s$); $(v,w)$ are coordinates in the $u$-dimensional
``unstable'' subspace (denoted $E^u$). The one-dimensional
coordinate $v$ corresponds to the eigenvalue $\lambda$ (and hence to
the one-dimensional ``direction of weakest expansion'' in $E^u$).

In the neighborhood $U$,
$$
\phi(t,(y,v,w))=(\exp(P_st)y;\exp(\lambda t)v,\exp(P_1t)w),
$$
and it follows from (\ref{l09}) that
\begin{equation}
 \label{l10}
|\exp(P_1t)w|\geq K\exp(\mu t)|w|,\quad t\geq 0.
\end{equation}
Denote by $E^u_1$ the one-dimensional invariant subspace
corresponding to $\lambda$.

We naturally identify $E^s\cap U$ and $E^u\cap U$ with the
intersections of $U$ with the corresponding local stable and
unstable manifolds of $p$, respectively.

Let us construct a special transverse section for the flow $\phi$.
We may assume that the point $r$ of nontransverse intersection of
$W^u(q)$ and $W^s(p)$ belongs to $U$. Take a hyperplane $\Sigma'$ in
$E^s$ of dimension $s-1$ that is transverse to the vector $X(r)$.
Set $\Sigma=\Sigma'+E^u$; clearly, $\Sigma$ is transverse to $X(r)$.

By a perturbation of the field $X$ outside $U$, we may get the
following: in a neighborhood of $r$, the component of intersection
$W^u(q)\cap \Sigma$ containing $r$ (for the perturbed field) has the
form of an affine space $r+L$, where $L$ is the tangent space,
$L=T_r(W^u(q)\cap\Sigma)$, of the intersection $W^u(q)\cap\Sigma$ at
the point $r$ for the unperturbed field (compare, for example, with
\cite{LeeSak}).

Let $\Sigma_r$ be a small transverse disk in $\Sigma$ containing the
point $r$. Denote by $\gam$ the component of intersection of
$W^u(q)\cap\Sigma_r$ containing $r$.
\medskip

{\bf Lemma 2. } {\em There exists $\ep>0$ such that if
$x\in\Sigma_r$ and}
\begin{equation}
\label{back} \mbox{dist}(\phi(t,x),O^-(r,\phi))<\ep,\quad t\leq 0,
\end{equation}
{\em then} $x\in\gam$.
\medskip

{\em Proof. } To simplify presentation, let us assume that $q$ is a
rest point; the case of a closed trajectory is considered using a
similar reasoning.

By the Grobman-Hartman theorem, there exists $\ep_0>0$ such that the
flow of $X$ in $B(2\ep_0,q)$ is topologically conjugate to the flow
of a linear vector field.

Denote by $A$ the intersection of the local stable manifold of $q$,
$W^s_{loc}(q)$, with the boundary of the ball $B(2\ep_0,q)$.

Take a negative time $T$ such that if $s=\phi(T,r)$, then
\begin{equation}
\label{back1} \phi(t,s)\in B(\ep_0,q),\quad t\leq 0.
\end{equation}
Clearly, if $\ep_0$ is small enough, then the compact sets $A$ and
$$
B=\{\phi(t,r):\; T\leq t\leq 0\}
$$
are disjoint. There exists a positive number $\ep_1<\ep_0$ such that
the $\ep_1$-neighborhoods of the sets $A$ and $B$ are disjoint as
well.

Take $\ep_2\in(0,\ep_1)$. There exists a neighborhood $V$ of the
point $s$ with the following property: if $y\in V\setminus
W^u_{loc}(q)$, then the first point of intersection of the negative
semitrajectory of $y$ with the boundary of $B(2\ep_0,q)$ belongs to
the $\ep_2$-neighborhood of the set $A$ (this statement is obvious
for a neighborhood of a saddle rest point of a linear vector field;
by the Grobman-Hartman theorem, it holds for $X$ as well).

Clearly, there exists a small transverse disk $\Sigma_s$ containing
$s$ and such that if $y\in \Sigma_s\cap W^u_{loc}(q)$, then the
first point of intersection of the positive semitrajectory of $y$
with the disk $\Sigma_r$ belongs to $\gam$ (in addition, we assume
that $\Sigma_s$ belongs to the chosen neighborhood $V$).

There exists $\ep\in(0,\ep_1-\ep_2)$ such that the flow of $X$
generates a local Poincar\'e transformation
$$
\sigma:\;\Sigma_r\cap B(\ep,r)\to\Sigma_s.
$$
Let us show that this $\ep$ has the desired property. It follows
from our choice of $\Sigma_s$ and (\ref{back}) with $t=0$ that if
$x\notin\gam$, then
$$
y:=\sigma(x)\in \Sigma_s\setminus W^u_{loc}(q);
$$
in this case, there exists $\tau<0$ such that the point
$z=\phi(\tau,y)$ belongs to the intersection of $B(\ep_2,A)$ with
the boundary of $B(2\ep_0,q)$. By (\ref{back1}),
\begin{equation}
\label{back2} \mbox{dist}(z,\phi(t,s))>\ep_0,\quad t\leq 0.
\end{equation}
At the same time,
\begin{equation}
\label{back3} \mbox{dist}(z,\phi(t,r))>\ep_1-\ep_2,\quad T\leq t\leq
0.
\end{equation}
Inequalities (\ref{back2}) and (\ref{back3}) contradict condition
(\ref{back}). Our lemma is proved.
\medskip

Now let us formulate the property of nontransversality of $W^u(q)$
and $W^s(p)$ at the point $r$ in terms of the introduced objects.

Let $\Pi^u$ be the projection to $E^u$ parallel to $E^s$.

The transversality of $W^u(q)$ and $W^s(p)$ at $r$ means that
$$
T_rW^u(q)+T_rW^s(p)=\RR^n.
$$
Since $\Sigma$ is a transverse section to the flow $\phi$ at $r$,
the above equality is equivalent to the equality
$$
L+E^s=\RR^n.
$$
Thus, the nontransversality means that
$$
L+E^s\neq \RR^n,
$$
which implies that
\begin{equation}
 \label{l11}
L':=\Pi^uL\neq E^u.
\end{equation}

We claim that there exists a linear isomorphism $J$ of $\Sigma$ for
which the norm $\|J-\mbox{Id}\|$ is arbitrarily small and such that
\begin{equation}
\label{iso} \Pi^uJL\cap E^u_1=\{0\}.
\end{equation}
Let $e$ be a unit vector of the line $E^u_1$. If $e\notin L'$, we
have nothing to prove (take $J=\mbox{Id}$). Thus, we assume that
$e\in L'$. Since $L'\neq E^u$, there exists a vector $v\in
E^u\setminus L'$.

Fix a natural number $N$ and consider a unit vector $v_N$ that is
parallel to $Ne+v$. Clearly, $v_N\to e$ as $N\to\infty$. There
exists a sequence $T_N$ of linear isomorphisms of $E^u$ such that
$T_Nv_N=e$ and
$$
\|T_N-\mbox{Id}\|\to 0,\quad N\to\infty.
$$
Note that $T^{-1}_Ne$ is parallel to $v_N$; hence, $T^{-1}_Ne$ does
not belong to $L'$, and
\begin{equation}
\label{iso1} T_N\Pi^uL\cap E^u_1=\{0\}.
\end{equation}
Define an isomorphism $J_N$ of $\Sigma$ by
$$
J_N(y,z)=(y,T_Nz)
$$
and note that
$$
\|J_N-\mbox{Id}\|\to 0,\quad N\to\infty.
$$
Let $L_N=J_N L$. Equality (\ref{iso1}) implies that
\begin{equation}
\label{iso2} \Pi^uL_N\cap E^u_1=\{0\}.
\end{equation}
Our claim is proved.

First we consider the case where $\mbox{dim}E^u\geq 2$. Since
$\mbox{dim}L'<\mbox{dim}E^u$ by (\ref{l11}) and $\mbox{dim}E^u_1=1$,
our reasoning above (combined with a Pugh-Robinson perturbation)
shows that we may assume that
\begin{equation}
\label{iso3} L'\cap E^u_1=\{0\}.
\end{equation}
For this purpose, we take a small transverse section $\Sigma'$
containing the point $r'=\phi(-1,r)$, denote by $\gam$ the component
of intersection of $W^u(q)$ with $\Sigma'$ containing $r'$, and note
that the local Poincar\'e transformation $\sigma$ generated by
$\Sigma'$ and $\Sigma$ takes $\gam$ to the linear space $L$ (in
local coordinates of $\Sigma$). The mapping $\sigma_N=J_N\sigma$ is
$C^1$-close to $\sigma$ for large $N$ and takes $\gam$ to $L_N$ for
which equality (\ref{iso2}) is valid. Thus, we get equality
(\ref{iso3}) for the perturbed vector field.

This equality implies that there exists a constant $C>0$ such that
if $(y;v,w)\in r+L$, then
\begin{equation}
 \label{l113}
|v|\leq C|w|.
\end{equation}

Fix $a>0$ such that $B(4a,p)\subset U$. Take a point
$\alpha=(0;a,0)\in E^u_1$ and a positive number $T$ and set
$\alpha_T=(r_y;a\exp(-\lambda T),0)$, where $r_y$ is the
$y$-coordinate of $r$. Construct a pseudotrajectory as follows:
$$
g(t)=\begin{cases}
\phi(t,r),\quad t\leq 0,\\
\phi(t,\alpha_T),\quad t>0.\\
\end{cases}
$$
Since
$$
|r-\alpha_T|=a\exp(-\lambda T)\to 0
$$
as $T \to \infty$, for any $d$ there exists $T$ such that $g$ is a
$d$-pseudotrajectory.
\medskip

{\bf Lemma 3. }{\em Assume that $b\in(0,a)$ satisfies the
inequality}
$$
\log K-\log
C+\left(\frac{\mu}{\lambda}-1\right)\left(\log\frac{a}{2} -\log
b\right)\geq 0.
$$
{\em Then for any $T>0$, reparametrization $h$, and a point $s\in
r+L$ such that $|r-s|<b$ there exists $\tau\in[0,T]$ such that}
$$
|\phi(h(\tau),s)-g(\tau)|\geq\frac{a}{2}.
$$
{\em Proof.} To get a contradiction, assume that
\begin{equation}
 \label{l13}
|\phi(h(\tau),s)-g(\tau)|<\frac{a}{2},\quad \tau\in[0,T].
\end{equation}
Let $s=(y_0;v_0,w_0)\in r+L$. Since $|r-s|<b$,
\begin{equation}
 \label{l14}
|v_0|<b.
\end{equation}
By (\ref{l13}),
$$
\phi(h(\tau),s)\in U,\quad \tau\in[0,T].
$$
Take $\tau=T$ in (\ref{l13}) to show that
$$
|v_0|\exp(\lambda h(T))>\frac{a}{2}.
$$
It follows that
\begin{equation}
 \label{l15}
h(T)>\lambda^{-1}\left(\log\frac{a}{2}-\log|v_0|\right).
\end{equation}
Set $\theta(\tau)=|\exp(P_1h(\tau))w_0|$; then $\theta(0)=|w_0|$. By
(\ref{l113}),
\begin{equation}
 \label{l16}
|v_0|\leq C\theta(0).
\end{equation}
By (\ref{l10}),
\begin{equation}
 \label{l17}
\theta(T)\geq K\exp(\mu h(T))\theta(0).
\end{equation}
We deduce from (\ref{l14})-(\ref{l17}) that
$$
\log\left(\frac{2\theta(T)}{a}\right)\geq \log
\theta(T)-\log|v_0\exp(\lambda h(T))|\geq
$$
$$
\geq \log K+\log \theta(0)-\log |v_0|+(\mu-\lambda)h(T)\geq
$$
$$
\geq \log K-\log C+\left(\frac{\mu}{\lambda}-1\right)
\left(\frac{a}{2}-\log|v_0|\right)\geq
$$
$$
\geq \log K-\log C+\left(\frac{\mu}{\lambda}-1\right)
\left(\frac{a}{2}-\log b\right)\geq 0.
$$
We get a contradiction with (\ref{l13}) for $\tau=T$ since the norm
of the $w$-coordinate of $\phi(h(T),s)$ equals $\theta(T)$, while
the $w$-coordinate of $g(T)$ is $0$. The lemma is proved.
\medskip

Let us complete the proof of Theorem 1 in case (B1). Assume that
$l,\delta_0>0$ are chosen for $\Sigma$ so that the LE holds.

Take $\epsilon\in(0,\min(\delta_0,\ep_0,a/2))$ so small that if
$|y-r|<\epsilon$, then $\phi(t,y)$ intersects $\Sigma$ at a point
$s$ such that
\begin{equation}
 \label{l011}
\mbox{dist}(\phi(t,s),r)<\ep_0,\quad |t|\leq l\ep.
\end{equation}
Consider the corresponding $d$ and a $d$-pseudotrajectory $g$
described above.

Assume that
\begin{equation}
 \label{l18}
\mbox{dist}(\phi(h(t),x),g(t))<\epsilon,\quad t\in \RR,
\end{equation}
for some point $x$ and reparametrization $h$ and set
$y=\phi(h(0),x)$.

Then $|y-r|<\ep$, and there exists a point $s=\phi(\tau,y)\in\Sigma$
with $|\tau|<l\ep$.

If $-l\ep\leq t\leq 0$, then
$$
\mbox{dist}(\phi(t,s),O^-(r,\phi))\leq\epsilon_0
$$
by (\ref{l011}).

If $t<-l\ep$, then $h(0)+\tau+t<h(0)$, and there exists $t'<0$ such
that $h(t')=h(0)+\tau+t$. In this case,
$$
\phi(t,s)=\phi(h(0)+\tau+t,x)=\phi(h(t'),x),
$$
and
$$
\mbox{dist}(\phi(t,s),O^-(r,\phi))\leq
\mbox{dist}(\phi(h(t'),x),\phi(t',r))\leq \epsilon_0.
$$

By Lemma 2, $s\in r+L$.  If $\epsilon$ is small enough, then
$|s-r|<b$, where $b$ satisfies the condition of Lemma 3, whose
conclusion contradicts (\ref{l18}).

This completes the consideration of case (B1) for
$\mbox{dim}W^u(p)\geq 2$. If $\mbox{dim}W^u(p)=1$, then the
nontransversality of $W^u(q)$ and $W^s(p)$ implies that $L\subset
E^s$. This case is trivial since any shadowing trajectory passing
close to $r$ must belong to the intersection $W^u(q)\cap W^s(p)$,
while we can contruct a pseudotrajectory ``going away" from $p$
along $W^u(p)$. If $\mbox{dim}W^u(p)=0$, $W^u(q)$ and $W^s(p)$
cannot have a point of nontransverse intersection.
\medskip

Case (B2). Passing from the vector field $X$ to $-X$, if necessary,
we may assume that $p$ is a closed trajectory. We ``linearize'' $X$
in a neighborhood of $p$ as described in the proof of Lemma 1 so
that the local Poincar\'e transformation of transverse section
$\Sigma$ is a linear mapping generated by a matrix $P$ with the
following properties: With respect to some coordinates in $\Sigma$,
\begin{equation}
 \label{l19}
P=\mbox{diag}(P_s,P_u),
\end{equation}
where $|\lambda_j|<1$ for the eigenvalues $\lambda_j$ of the matrix
$P_s$, and $|\lambda_j|>1$ for the eigenvalues $\lambda_j$ of the
matrix $P_u$, every eigenvalue has multiplicity 1, and $P$ is in a
Jordan form.

The same reasoning as in case (B1) shows that it is possible to
perform such a ``linearization" (and other perturbations of $X$
performed below) so that the nontransversality of $W^u(q)$ and
$W^s(p)$ is preserved.

Consider an eigenvalue $\lambda$ of $P_u$ such that
$|\lambda|\leq|\mu|$ for the remaining eigenvalues $\mu$ of $P_u$.

We treat separately the following two cases.

Case (B2.1): $\lambda\in \RR$.

Case (B2.2): $\lambda\in \CC\setminus \RR$.
\medskip

Case (B2.1). Applying a perturbation, we may assume that
$$
P_u=\mbox{diag}(\lambda,P_1),
$$
where $|\lambda|<|\mu|$ for the eigenvalues $\mu$ of the matrix
$P_1$ (thus, there exists a one-dimensional direction of ``weakest
expansion'' in $W^u(p)$). In this case, we apply precisely the same
reasoning as that applied to treat case (B1) (we leave details to
the reader).
\medskip

Case (B2.2). Applying one more perturbation of $X$, we may assume
that
$$
\lambda=\nu+i\eta=\rho\exp\left(\frac{2\pi m_1i}{m}\right),
$$
where $m_1$ and $m$ are relatively prime natural numbers, and
$$
P_u=\mbox{diag}(Q,P_1),
$$
where
$$
Q=\left(
\begin{matrix}
\nu&-\eta\\
\eta&\nu\\
\end{matrix}
\right)
$$
with respect to some coordinates $(y,v,w)$ in $\Sigma$, where
$\rho=|\lambda|<|\mu|$ for the eigenvalues $\mu$ of the matrix
$P_1$.

Denote
$$
E^s=\{(y,0,0)\},\quad E^u=\{(0,v,w)\},\quad E^u_1=\{(0,v,0)\}.
$$
Thus, $E^s$ is the ``stable subspace," $E^u$ is the ``unstable
subspace," and $E^u_1$ is the two-dimensional ``unstable subspace of
the weakest expansion."

Geometrically, the Poincar\'e transformation
$\sigma:\Sigma\to\Sigma$ (extended as a linear mapping to $E^u_1$)
acts on $E^u_1$ as follows: the radius of a point is multiplied by
$\rho$, while $2\pi m_1/m$ is added to the polar angle.

As in the proof of Lemma 1, we take a small neighborhood $W$ of the
origin of the transverse section $\Sigma$ so that, for points $x\in
W$, the function $\al(x)$ (the time of first return to $\Sigma$) is
defined.

We assume that the point $r$ of nontransverse intersection of
$W^u(q)$ and $W^s(p)$ belongs to the section $\Sigma$. Similarly to
case (B1), we perturb $X$ so that, in a neighborhood of $r$, the
component of intersection of $W^u(q)\cap\Sigma$ containing $r$ has
the form of an affine space, $r+L$.

Let $\Pi^u$ be the projection in $\Sigma$ to $E^u$ parallel to
$E^s$, and let $\Pi^u_1$ be the projection to $E^u_1$; thus,
$$
\Pi^u(y,u,v)=(0,u,v)\mbox{ and }\Pi^u_1(y,u,v)=(0,u,0).
$$

The nontransversality of $W^u(q)$ and $W^s(p)$ at $r$ means that
$$
L'=\Pi^u L\neq E^u
$$
(see case (B1)). Applying a reasoning similar to that in case (B1),
we perturb $X$ so that if $L''=L'\cap E^u_1$, then
$$
\mbox{dim}L''<\mbox{dim}E^u_1=2.
$$
Hence, either $\mbox{dim}L''=1$ or $\mbox{dim}L''=0$. We consider
only the first case, the second one is trivial.

Denote by $A$ the line $L''$. Images of $A$ under degrees of
$\sigma$ (extended to the whole plane $E^u_1$) are $m$ different
lines in $E^u_1$.

In what follows, we refer to an obvious geometric statement (given
without a proof).
\medskip

{\bf Proposition 1.} {\em Consider Euclidean space $\RR^n$ with
coordinates $(x_1,\dots,x_n)$. Let $x'=(x_1,x_2)$,
$x''=(x_3,\dots,x_n)$, and let $G$ be the plane of coordinate $x'$.
Let $D$ be a hyperplane in $\RR^n$ such that
$$
D\cap G=\{x_2=0\}.
$$
For any $b>0$ there exists $c>0$ such that if $x=(x',x'')\in D$ and
$x'=(x'_1,x'_2)$, then either $|x'_2|\leq b|x'_1|$ or} $|x''|\geq
c|x'|$.
\medskip

Take $a>0$ such that the $2a$-neighborhood of the origin in $\Sigma$
belongs to $W$. We may assume that if $v=(v_1,v_2)$, then the line
$A$ is $\{v_2=0\}$.

Take $b>0$ such that the images of the cone
$$
C=\{v:\;|v_2|\leq b|v_1|\}
$$
in $E^u_1$ under degrees of $\si$ intersect only at the origin
(denote these images by $C_1,\dots,C_m)$.

We apply Proposition 1 to find a number $c>0$ such that if
$(0,v,w)\in L'$, then either $(0,v,0)\in C$ or
\begin{equation}
\label{12-13} |w|\geq c|v|.
\end{equation}

Take a point $\beta=(0,v,0)\in\Sigma$, where $|v|=a$, such that
$\beta\notin C_1\cup\dots\cup C_m$.

For a natural number $N$, set $\beta_N=(r_y,P_u^{-N}(v,0))\in\Sigma$
(we recall that equality (\ref{l19}) holds), where $r_y$ is the
$y$-coordinate of $r$. We naturally identify $\beta$ and $\beta_N$
with points of $M$ and consider the following pseudotrajectory:
$$
g(t)=\begin{cases}
\phi(t,r),\quad t\leq 0;\\
\phi(t,\beta_N),\quad t>0.\\
\end{cases}
$$

The following statement (similar to Lemma 2) holds: there exists
$\epsilon_0>0$ such that if
$$
\mbox{dist}(\phi(t,s),O^-(r,\phi))<\epsilon_0,\quad t\leq 0,
$$
for some point $s\in \Sigma$, then $s\in r+L$.

Since $\beta$ does not belong to the closed set $C_1\cup\dots\cup
C_m$, we may assume that the disk in $E^u_1$ centered at $\beta$ and
having radius $\epsilon_0$ does not intersect the set
$C_1\cup\dots\cup C_m$.

Define numbers
\begin{multline*}
\al_1(N)=\al(\beta_N),\;\al_2(N)=\al_1(N)+\al(\sigma(\beta_N)),
\dots,\\ \al_N(N)=\al_{N-1}(N)+\al(\sigma^{N-1}(\beta_N)).
\end{multline*}

Take $\de_0$ and $l$ for which LE holds for the neighborhood $W$
(reducing $W$, if necessary). Take
$\epsilon<\min(\epsilon_0/l,\delta_0)$ and assume that there exists
the corresponding $d$ (from the definition of the class OrientSh).
Take $N$ so large that $g$ is a $d$-pseudotrajectory.

Let $h$ be a reparametrization; assume that
$$
|\phi(h(t),p_0)-g(t)|<\ep,\quad 0\leq t\leq \al_N(N),
$$
for some point $p_0\in\Sigma$.

Since $g(\al_k(N))\in \Sigma$ for $0\leq k\leq N$ by construction,
there exist numbers $\chi_k$ such that
$$
\left|\sigma^{\chi_k}(p_0)-g(\al_k(N))\right|<\ep_0,\quad 0\leq
k\leq N.
$$

To complete the proof of Theorem 1, let us show that for any $p_0\in
r+L$ and any reparametrization $h$ there exists $t\in[0,\al_N(N)]$
such that
$$
\mbox{dist}(\phi(h(t),p_0),g(t))\geq\epsilon.
$$
Assuming the contrary, we see that
$$
\left|\sigma^{\chi_k}(p_0)-g(\al_k(N))\right|<\epsilon_0,\quad 0\leq
k\leq N,
$$
where the numbers $\chi_k$ were defined above.

We consider two possible cases.

If
$$
\Pi^u_1p_0\in C
$$
($C$ is the cone defined before estimate (\ref{12-13})), then
$$
\Pi^u_1\si^{\chi_k}(p_0)\in C_1\cup\dots\cup C_m.
$$
By construction, $\Pi^u_1g(\al_N(N))$ is $\beta$. Hence,
$$
\left|\Pi^u_1\sigma^{\chi_N}(p_0)-\Pi^u_1g(\al_N(N))\right|>\epsilon_0,
$$
and we get the desired contradiction.

If
$$
\Pi^u_1p_0 \notin C
$$
and $p_0=(y_0,v_0,w_0)$, then $(0,v_0,w_0)\in L'$, and it follows
from (\ref{12-13})) that $|w_0|\geq c|v_0|$. In this case,
decreasing $\ep_0$, if necessary, we apply the reasoning similar to
Lemma 3.

Thus, we have shown that
\begin{equation}
\label{aaa} \Int^1(\OrientSh\setminus\Bb)\subset\Int^1(\KS)=\Ss.
\end{equation}

It was shown in \cite{PilFlow} that $\Ss\subset\StSh$; since the set
$\Ss$ is $\Cone$-open and $\Ss \cap \Bb = \emptyset$,
\begin{equation}
\label{aab} \Ss\subset\Int^1(\StSh\setminus\Bb)\subset
\Int^1(\OrientSh\setminus\Bb).
\end{equation}
Inclusions (\ref{aaa}) and (\ref{aab}) prove Theorem 1.
\medskip

By Lemma 1, if $X\in \Int^1(\OrbitSh)$, then $X\in \Int^1(\Tt)$. For
nonsingular flows, the latter inclusion implies that $X$ is
$\Omega$-stable \cite{GanWen} (note that this is not the case for
flows with rest points \cite{Mane}). Now, based on the second part
of the proof of Theorem 1, one easily proves Theorem 3 following the
same lines as in \cite[Theorem 4]{PilRodSak}.

\section{Proof of Theorem 2}

Consider a vector field $X^*$ on the manifold $M = S^2 \times S^2$
that has the following properties (F1)-(F3) ($\phi^*$ denotes the
flow generated by $X^*$).

\begin{enumerate}
\item[] \text{(F1)} The nonwandering set of $\phi^*$ is the union of
four rest points $p^*, q^*, s^*, u^*$.
\item[] \text{(F2)} For some
$\delta > 0$ we can introduce coordinates in the neighborhoods
$B(\delta, p^*)$ and $B(\delta, q^*)$ such that
$$
X^*(x)=J^*_p(x-p^*),\quad x \in B(\delta, p^*),\quad \mbox{and}
\quad X^*(x)=J^*_q(x-q^*),\quad x \in B(\delta, q^*),
$$
where
\begin{equation}\notag
J^*_p = -J^*_q = \left(
            \begin{array}{cccc}
              -1 & 0 & 0 & 0 \\
              0 & -2 & 0 & 0 \\
              0 & 0 & 1 & -1 \\
              0 & 0 & 1 & 1 \\
            \end{array}
          \right),
\end{equation}


\item[] \text{(F3)} The point $s^*$ is an attracting hyperbolic rest point. The point
$u^*$ is a repelling hyperbolic rest point. The following condition
holds:
\begin{equation}\label{Text16.5}
\Wu(p^*) \setminus \{p^*\} \subset \Ws(s^*), \quad \Ws(q^*)
\setminus \{q^*\} \subset \Wu(u^*).
\end{equation}

The intersection of $\Ws(p^*) \cap \Wu(q^*)$ consists of a single
trajectory $\al^*$, and for any $x \in \al^*$, the condition
\begin{equation}\label{dim3}
\dim T_x \Ws(p^*) \oplus T_x\Wu(q^*) = 3
\end{equation}
holds.
\end{enumerate}

These conditions imply that the two-dimensional manifolds $\Ws(p^*)$
and $\Wu(q^*)$ intersect along a one-dimensional curve in the
four-dimensional manifold $M$. Thus, $\Ws(p^*)$ and $\Wu(q^*)$ are
not transverse; hence, $X^* \in \Bb$.

A construction of such a vector field is given in the Appendix.

To prove Theorem 2, we show that $X^*\in\Int^1(\OrientSh)$.
\medskip

The vector field $X^*$ satisfies Axiom A and the no-cycle condition;
hence, $X^*$ is $\Omega$-stable. Thus, there exists a neighborhood
$V$ of $X^*$ in the $C^1$-topology such that for any field $X \in
V$, its nonwandering set consists of four hyperbolic rest points $p,
q, s, u$ which belong to small neighborhoods of $p^*, q^*, s^*,
u^*$, respectively. We denote by $\phi$ the flow of any $X \in V$
and by $W^s(p),W^u(p)$ etc the corresponding stable and unstable
manifolds.

Note that if the neighborhood $V$ is small enough, then there exists
a number $c>0$ (the same for all $X \in V$) such that
$$
B(c, s^*) \subset W^s(s)\quad\mbox{and} \quad B(c, u^*) \subset
W^u(u).
$$
Consider the set $\Theta = \Wu(p^*) \cap \partial B(\delta, p^*)$
(where $\partial A$ is the boundary of a set $A$). Condition
\sref{Text16.5} implies that there exists a neighborhood $U_\Theta$
of $\Theta$ and a number $T>0$ such that
$$
\phi^*(T, x) \in B(c/2, s^*),\quad x \in U_\Theta.
$$
Reducing $V$, if necessary, we may assume that
$$
\Wu(p) \cap \partial B(\delta, p) \subset U_\Theta\quad\mbox{and}
\quad \phi(T, x) \in B(c, s^*), \quad x \in U_\Theta.
$$
Hence, $\Wu(p)\setminus\{p\} \subset \Ws(s)$, and
\begin{equation}\label{Text3.5}
\Wu(p) \cap \Ws(q) = \emptyset.
\end{equation}
Similarly, we may assume that $\Ws(q)\setminus\{q\}\subset \Wu(u)$.

The following two cases are possible for $X \in V$.
\begin{enumerate}
\item[] \text{(S1)} $\Ws(p) \cap \Wu(q) = \emptyset$.
\item[] \text{(S2)} $\Ws(p) \cap \Wu(q) \ne \emptyset$.
\end{enumerate}

In case \text{(S1)}, $X$ is a Morse-Smale field; hence, $X\in \Ss$.
Since $\Ss\subset\StSh$ (see \cite{PilFlow}), $X\in\OrientSh$.
\medskip

{\bf Remark 1. } In fact, it is shown in \cite{PilFlow} that if a
vector field $X\in \Ss$ does not have closed trajectories (as in our
case), then $X$ has the Lipschitz shadowing property without
reparametrization of shadowing trajectories: there exists $L>0$ such
that if $g(t)$ is a $d$-pseudotrajectory with small $d$, then there
exists a point $x$ such that
$$
\dist(g(t),\phi(t,x))\leq Ld,\quad t\in\RR.
$$
We refer to this fact below.
\medskip

Thus, in the rest of the proof of Theorem 2, we consider case
\text{(S2)}. Our goal is to show that if the neighborhood $V$ is
small enough, then $X\in\OrientSh$.
\medskip

{\bf Lemma 4. }{\em If the neighborhood $V$ is small enough, then
the intersection $\Ws(p) \cap \Wu(q)$ consists of a single
trajectory}.
\medskip

{\em Proof. } Denote $x_p^* = \al^* \cap \partial B(\delta, p^*)$
and $x_q^* = \al^* \cap \partial B(\delta, q^*)$.

Consider  sections $Q_p$ and $Q_q$ transverse to $\al$ at the points
$x_p^*$ and $x_q^*$, respectively,  and the corresponding Poincar\'e
map $F^*: Q_q \to Q_p$. Consider the curves $\xi_p^* = \Ws(p^*) \cap
Q_p \cap B(\delta/2, x_p^*)$ and $\xi_q^* = \Ws(q^*) \cap Q_q \cap
B(\delta/2, x_q^*)$. Note that $\xi^*_p$ and $F^*(\xi^*_q)$
intersect at a single point $x^*_p$.

Let $\xi_p = \Ws(p) \cap Q_p \cap B(\delta/2, x^*_p)$ and $\xi_q =
\Wu(q) \cap Q_q \cap B(\delta/2, x^*_q)$. Let $F$ be the Poincar\'e
transformation for $X$ from $Q_q$ to $Q_p$ similar to $F^*$.

If the neighborhood $V$ is small enough, then the curves $\xi_p$,
$\xi_q$, and $F(\xi_q)$ are $\Cone$-close to $\xi^*_p$, $\xi^*_q$,
and $F^*(\xi^*_q)$, respectively (hence, the intersection of $\xi_p$
and $F(\xi_q)$ contains not more than one point).

The same reasoning as in the proof of (\ref{Text3.5}) shows that if
the neighborhood $V$ is small enough, $x\in\Ws(p)\setminus\{p\}$,
and the trajectory of $x$ does not intersect $\xi_p$, then
$x\in\Wu(u)$.

Thus, any trajectory in $\Ws(p) \cap \Wu(q)$ must intersect $\xi_p$;
similarly, it must intersect $\xi_q$ as well as $F(\xi_q)$.

It follows that the intersection $\Ws(p) \cap \Wu(q)$ (which is
nonempty since we consider case \text{(S2)}) consists of a single
trajectory containing the unique point $x_p$ of intersection of
$\xi_p$ and $F(\xi_q)$ (we denote this trajectory by $\alpha$). This
completes the proof of Lemma 4.
\medskip

{\bf Remark 2. } Let us note an important property of intersection
of $\Ws(p)$ and $\Wu(q)$ along $\alpha$ (see (\ref{Add3.2}) below).

Let $x_q=F^{-1}(x_p)$; denote by $i_p$ and $i_q$ unit tangent
vectors to the curves $\xi_p$ and $\xi_q$ at $x_p$ and $x_q$,
respectively. Our reasoning above and condition \sref{dim3} show
that if the neighborhood $V$ is small enough, then the vectors $i_p$
and $D F(x_q) i_q$ are not parallel:
\begin{equation}\label{Add3.2'}
D F(x_q) i_q \nparallel i_p.
\end{equation}

Take any two points $y_p=\phi(t_1,x_p)$ and $y_q=\phi(t_2,x_q)$ with
$t_1\geq0,t_2\leq 0$; let $S_p$ and $S_q$ be smooth transversals to
$\al$ at these points. Let $e_p$ and $e_q$ be tangent vectors of
$S_p \cap \Ws(p)$ and $S_q \cap \Wu(q)$ at $y_p$ and $y_q$,
respectively. Denote by $f: S_q \to S_p$, $H_p: Q_p \to S_p$, and
$H_q: S_q \to Q_q$ the corresponding Poincar\'e transformations for
$X$. Then $f = H_p \circ F \circ H_q$,
$$
e_p \parallel D H_p(x_p)i_p, \quad\mbox{and}\quad e_q \parallel D
H_q^{-1}(x_q) i_q.
$$
Hence, $D f (y_q) e_q
\parallel D H_p \circ D F(x_q) i_q$, and it follows from
\sref{Add3.2'} that
\begin{equation}\label{Add3.2}
D f(y_q) e_q \nparallel e_p.
\end{equation}

Now it remains to show that if $V$ is small enough and $X\in V$,
then $X\in\OrientSh$ (recall that we consider case (S2)). This proof
is rather complicated, and we first describe its scheme.

We fix two points $y_p,y_q\in\alpha$ in small neighborhoods $U_p$
and $U_q$ of $p$ and $q$, respectively (the choice of $U_p$ and
$U_q$ is specified later). We consider special pseudotrajectories
(of type \text{Ps}): the "middle"{} part of such a pseudotrajectory
is the part of $\alpha$ between $y_q$ and $y_p$, while its
"negative"{} and "positive"{} tails are parts of trajectories that
start near $y_q$ and $y_p$, respectively. We show that our shadowing
problem is reduced to shadowing of pseudotrajectories of type
\text{Ps}.

The key part of the proof is a statement "on four balls."{} It is
shown that if $B_1,\dots,B_4$ are small balls such that $B_1$ and
$B_4$ are centered at points of $W^s(q)$ and $W^u(p)$, while $B_2$
and $B_3$ are centered at $y_q$ and $y_p$, respectively, then there
exists an exact trajectory that intersects $B_1,\dots,B_4$
successfully as time grows. This statement (and its analog) allows
us to prove that pseudotrajectories of type \text{Ps} can be
shadowed.
\medskip

Let us fix points $y_p, y_q \in \al$ (everywhere below, we assume
that $y_p=\alpha(T_p)$ and $y_q=\alpha(T_q)$ with $T_p > T_q$) and a
number $\delta> 0$. We say that ${g}(t)$ is a pseudotrajectory of
type \text{Ps}$(\delta)$ if
\begin{equation}\label{Pst1.2.2}
{g}(t) =
\begin{cases}
\phi(t - T_p, x_p), & t > T_p, \\
\phi(t - T_q, x_q), & t < T_q, \\
\alpha(t), & t \in [T_q, T_p],
\end{cases}
\end{equation}
for some points
$$
 x_p \in B(\delta, y_p)\quad \mbox{and} \quad
x_q \in B(\delta, y_q).
$$

Fix an arbitrary $\ep > 0$. We prove the following two statements
(Propositions 2 and 3). In these statements, we say that a
pseudotrajectory $g(t)$ can be $\ep$-shadowed if there exists a
reparametrization $h$ and a point $p$ such that (\ref{02}) holds.

An $\Omega$-stable vector field has a continuous Lyapunov function
that strictly decreases along wandering trajectories (see
\cite{PughShub}). Hence, there exist small neighborhoods $U_p$ and
$U_q$ of points $p$ and $q$, respectively, such that
\begin{equation}\label{Add3.4}
\phi(t,x)\notin U_q,\quad x\in U_p,\;t\geq 0.
\end{equation}

\medskip

{\bf Proposition 2. } {\em For any $\delta > 0$, $y_p \in  \al\cap
U_p$, and $y_q \in \al\cap U_q$ there exists $d> 0$ such that if
$g(t)$ is a $d$-pseudotrajectory of $X$, then either $g(t)$ can  be
$\ep$-shadowed or there exists a pseudotrajectory $g^*(t)$ of type}
\text{Ps}$(\delta)$ {\em with these $y_p$ and $y_q$ such that}
$\dist(g(t), g^*(t)) < \ep/2,\quad t\in\RR$.
\medskip

{\bf Proposition 3. } {\em There exists $\delta > 0$, $y_p
\in\al\cap U_p$, and $y_q \in \al\cap U_q$ such that any
pseudotrajectory of type} \text{Ps}$(\delta)$ {\em with these $y_p$
and $y_q$ can be $\ep/2$-shadowed.}
\medskip

Clearly, Propositions 2 and 3 imply that $X \in \OrientSh$.
\medskip

To prove Proposition 2, we need an auxiliary statement.
\medskip

{\bf Lemma 5.} {\em For any $x \in \al$ and $\ep, \ep_1 > 0$ there
exists $d>0$ such that if
\begin{equation}\label{constr1.1}
\{g(t): \; t \in \RR\} \cap B(\ep_1, x) = \emptyset,
\end{equation}
for a $d$-pseudotrajectory $g(t)$, then one can find $x_0 \in M$ and
$h(t) \in \Rep$ such that}
$$
\dist(g(t), \phi(h(t), x_0)) < \ep, \quad t \in \RR.
$$
\medskip

{\em Proof. } Take $\Delta < \ep_1/2$ such that if $a_p=\phi(1,x)$
and $a_q=\phi(-1,x)$, then $a_p, a_q \notin  B(\Delta, x)$. Let
$S_p$ and $S_q$ be three-dimensional transversals to $\alpha$ at
$a_p$ and $a_q$, respectively. Let ${f : S_q \to S_p}$ be the
corresponding Poincar\'e mapping. Note that the intersections
$\Wu(q) \cap S_q$ and $\Ws(p) \cap S_p $ near $a_q$ and $a_p$ are
one-dimensional, hence the curves $f(\Wu(q) \cap S_q)$ and $\Ws(p)
\cap S_p$ in $S_p$ are nontransverse.

It is shown in \cite{MorSakSumi, PughRob} that there exists an
arbitrarily small perturbation of the field $X$ supported in
$B(\Delta, x)$ and such that the Poincar\'e mapping $\tilde{f}: S_q
\to S_p$ of the perturbed field $\tilde{X}$ satisfies the condition
$$
\tilde{f}(\Wu(q) \cap S_q) \cap (\Ws(p) \cap S_p) = \emptyset.
$$
Similarly to case (S1), we conclude that we can find $\tilde{X} \in
\Ss$.

Set $\ep_2 = \min(\ep, \ep_1/2)$ and find $d> 0$ such that any
$d$-pseudotrajectory of the field $\tilde{X}$ can be
$\ep_2$-shadowed. We assume, in addition, that
\begin{equation}
\label{lem1} \Delta + d < \ep_1.
\end{equation}
Consider an arbitrary $d$-pseudotrajectory $g(t)$ of $X$ for which
\sref{constr1.1} holds. By (\ref{lem1}), $g(t)$ is a
$d$-pseudotrajectory of the field $\tilde{X}$. Due to the choice of
$d$, there exists $x_0 \in M$ and $h(t) \in \Rep$ such that
$$
\dist(g(t), \tilde{\phi}(h(t), x_0)) < \ep_2,
$$
where $\tilde{\phi}$ is the flow of $\tilde{X}$. Hence,
$\{\tilde{\phi}(h(t), x_0), \; t \in \RR\} \cap B(\ep_1, x) =
\emptyset$; it follows that $\tilde{\phi}(h(t), x_0) = \phi(h(t),
x_0)$, which proves Lemma 5.
\medskip

{\em Proof of Proposition 2. } Take $\delta > 0$,  $y_p
\in\alpha\cap U_p$, and $y_q \in \alpha\cap U_q$. Let
$y_q=\alpha(T_q)$ and $y_p=\alpha(T_p)$. There exists
$\delta_1\in(0,\min(\delta,\ep))$ such that $B(\delta_1, y_p)
\subset U_p$, $B(\delta_1, y_q) \subset U_q$, and if $x_p\in
B(\delta_1, y_p)$ and $x_q\in B(\delta_1, y_q)$, then
\begin{equation}
\label{psu} g^*(t) =
\begin{cases}
\phi(t-T_p, x_p), & t > T_p, \\
 \al(t), & t \in [T_q, T_p], \\
 \phi(t-T_q, x_q), & t < T_q,
\end{cases}
\end{equation}
is a pseudotrajectory of type \text{Ps}$(\delta)$.

Take $x=\alpha(T)$, where $T\in(T_q,T_p)$. Applying Lemma 5, we can
find $\ep_1>0$ such that if $d$ is small enough, then for any
$d$-pseudotrajectory $g(t)$, one of the following two cases holds
(after a shift of time):
\begin{itemize}
\item[] \text{(A1)}
$$
\{g(t), \; t \in \RR\} \cap B(\ep_1, x) = \emptyset,
$$
and $g(t)$ can be $\ep$-shadowed;
\item[] \text{(A2)}
$$
g(T_p) \in B(\delta_1/2, y_p), \quad g(T_q) \in B(\delta_1/2, y_q),
$$
and
$$
\dist(g(t), \alpha(t)) < \ep/2,\quad t\in[T_q,T_p].
$$
\end{itemize}
To prove Proposition 2, it remains to consider case \text(A2).

Apply the same reasoning as in Lemma 5 to construct a field
$\tilde{X}\in \Ss$ that coincides with $X$ outside $B(\delta_1/2,
y_q)$; let $\tilde{\phi}$ be the flow of $\tilde{X}$.

Note that $\tilde{X}$ does not have closed trajectories. Reducing
$d$, if necessary, we may assume that any $d$-pseudotrajectory of
$\tilde{X}$ can be $\delta_1/2$-shadowed in the sense of Remark 1.

Consider the mapping
$$
\tilde{g}_p(t) =
\begin{cases}
\tilde{\phi}(t - T_p, g(T_p)), & t < T_p, \\
g(t), & t \in [T_p, T], \\
\tilde{\phi}(t - T, g(T)), & t > T,
\end{cases}
$$
where
$$
T = \inf\{t > T_p: \tilde{g}_p(t) \in B(\delta_1, y_q)\}
$$
(if $\{t> T_p: \tilde{g}_p(t) \in B(\delta_1, y_q)\} = \emptyset$,
we set $T = +\infty$). Since
$$
B(\delta_1/2,g(t))\cap B(\delta_1/2, y_q)=\emptyset
$$
for $t\in[T_p,T)$, $\tilde{g}_p(t)$ is a $d$-pseudotrajectory of
$\tilde{X}$. Hence, there exists a point $x_p$ such that
\begin{equation}\label{Text6.1}\notag
\dist(\tilde{g}_p(t), \tilde{\phi}(t - T_p, x_p)) < \delta_1/2,
\quad t \in \RR.
\end{equation}
The first inclusion in (A2) implies that $x_p \in B(\delta, y_p)$.

Since trajectories of $X$ and $\tilde{X}$ coincide outside
$B(\delta_1/2, y_q)$, we deduce from \sref{Add3.4} that $T =
+\infty$; hence,
$$
\dist(g(t), \phi(t - T_p, x_p)) < \delta_1/2, \quad t \geq T_p.
$$

Similarly (reducing $d$, if necessary), we find $x_q \in B(\delta,
y_q)$ such that
$$
\dist(g(t), \phi(t - T_q, x_q)) < \delta_1/2, \quad t \leq T_q.
$$

Clearly, the mapping (\ref{psu}) is a pseudotrajectory of type
\text{Ps}$(\delta)$ such that
$$
\dist(g(t), g^*(t)) < \ep/2,\quad t\in\RR.
$$
This completes the proof of Proposition 2.

In the remaining part of the paper, we prove Proposition 3. Let us
recall that we consider a vector field $X$ in a small neighborhood
$V$ of $X^*$ for which $\Ws(p)\cap\Wu(q)\neq\emptyset$.

Without loss of generality, we may assume that
$$
O^+(B(\ep/2, s),\phi) \subset B(\ep, s)\quad\mbox{and} \quad
O^-(B(\ep/2, u),\phi) \subset B(\ep, u).
$$

Take  $m \in (0, \ep/8)$ such that $B(m, p) \subset U_p$, $B(m, q)
\subset U_q$ and the flow of the vector field $X$ in the
neighborhoods $B(2m, p)$ and $B(2m, q)$ is conjugate by a
homeomorphism to the flow of a linear vector field.

We take points $y_p=\alpha(T_p) \in B(m/2, p)\cap \al$ and
$y_q=\alpha(T_q) \in B(m/2, q) \cap \al$. Then $O^+(y_p,\phi)
\subset B(m, p)$ and $O^-(y_q,\phi)\subset B(m, q)$. Take $\delta>0$
such that if $g(t)$ is a pseudotrajectory of type
\text{Ps}$(\delta)$ (with $y_p$ and $y_q$ fixed above), $t_0 \in
\RR$, and $x_0 \in B(2\delta, g(t_0))$, then
\begin{equation}\label{Pst1.4.1}
\dist(\phi(t-t_0, x_0), g(t)) < \ep/2 , \quad |t-t_0| \leq T+1,
\end{equation}
where $T=T_p - T_q$.

Consider a number $\tau > 0$ such that if $x \in \Wu(p) \setminus
B(m/2, p)$, then $\phi(\tau, x) \in B(\ep/8,s)$. Take
$\ep_1\in(0,m/4)$ such that if two points $z_1, z_2 \in M$ satisfy
the inequality $\dist(z_1, z_2) < \ep_1$, then
$$
\dist(\phi(t, z_1), \phi(t, z_2))< \ep/8, \quad |t| \leq \tau.
$$
In this case, for any $y \in B(\ep_1, x)$ (recall that we consider
$x \in \Wu(p)\setminus B(m/2, p)$), the following inequalities hold:
\begin{equation}\label{Pst1.5.1}
\dist(\phi(t, x), \phi(t, y)) < \ep/4, \quad t \geq 0.
\end{equation}
Reducing  $\ep_1$, if necessary, we may assume that if $x'\in \Ws(q)
\setminus B(m/2, q)$ and $y' \in B(\ep_1, x')$, then
\begin{equation}\label{Pst1.5.1'}\notag
\dist(\phi(t, x'), \phi(t, y')) < \ep/4, \quad t \leq 0.
\end{equation}

Let $g(t)$ be a pseudotrajectory of type \text{Ps}$(\delta)$, where
$\delta$, $y_p$, and $y_q$ satisfy the above-formulated conditions.
We claim that if $\delta$ is small enough, then $g(t)$ can be
$\ep/2$-shadowed (in fact, we have to reduce $\delta$ and to impose
additional conditions on $y_p$ and $y_q$). Below we denote
$W^u_{loc}(p, m)=\Wu(p) \cap B(m, p)$ etc.

Additionally decreasing $\delta$, we may assume that for any points
$z_p \in W^u_{loc}(p, m)$, $x_0 \in B(\delta, y_p)$, and $s>0$ such
that $\phi(s, x_0) \in B(\delta, z_p)$, the following inclusions
hold:
\begin{equation}\label{Text23.1}
\phi(t, x_0) \in B(2m, p), \quad t \in [0, s].
\end{equation}

Let us consider several possible cases.
\medskip

Case (P1):  $x_p \notin \Ws(p)$ and  $x_q \notin \Wu(q)$. Let
$$
T' = \inf \{ t \in \RR: \; \phi(t, x_p) \notin B(p, 3m/4) \}.
$$
If $\delta$ is small enough, then $\dist(\phi(T',x_p), \Wu(p)) <
\ep_1$. In this case, there exists a point $z_p \in W^u_{loc}(p,
m)\setminus B(m/2, p)$ such that
\begin{equation}\label{Text8.0.5}
\dist(\phi(T', x_p), z_p) < \ep_1.
\end{equation}

Applying a similar reasoning in a neighborhood of $q$ (and reducing
$\delta$, if necessary), we find a point $z_q \in W^s_{loc}(q,
m)\setminus B(m/2, q)$ and a number $T'' < 0$ such that
$\dist(\phi(T'', x_q), z_q) < \ep_1$.

Let us formulate a key lemma which we prove later (precisely this
lemma is the above-mentioned statement "on four balls").
\medskip

{\bf Lemma 6. }{\em There exists  $m >0$ such that for any points
$$
y_p \in B(m, p) \cap \al, \quad z_p \in W^u_{loc}(p, m) \setminus
\{p\},
$$
$$
y_q \in B(m, q) \cap \al, \quad z_q \in W^s_{loc}(q, m) \setminus
\{q\},
$$
and for any number $m_1 > 0$ there exists a trajectory of the vector
field $X$ that intersects successively the balls $ B(m_1, z_q)$,
$B(m_1, y_q)$, $B(m_1, y_p)$, and $B(m_1, z_p)$ as time grows.}
\medskip

We reduce $m$ to satisfy Lemma 6 and apply this lemma with $m_1 =
\min(\delta,\ep_1)$. Find a point $x_0$ and numbers $t_1< t_2 < t_3
< t_4$ such that
$$
\phi(t_1, x_0) \in B(m_1, z_q), \quad \phi(t_2, x_0) \in B(m_1,
y_q),
$$
$$
\phi(t_3, x_0) \in B(m_1, y_p), \quad \phi(t_4, x_0) \in B(m_1,
z_p).
$$
Inequalities \sref{Pst1.4.1} imply that if $\delta$ is small enough,
then
\begin{equation}\label{Pst1.6.5}
\dist(\phi(t_3 + t, x_0), g(T_p + t)) < \ep/2, \quad t \in [T_q -
T_p, 0].
\end{equation}
Define a reparametrization $h(t)$ as follows:
$$
h(t) = \begin{cases}
h(T_q + T'' + t) = t_1 + t, & t<0,\\
h(T_p + T' + t) = t_4 + t, & t>0, \\
h(T_p + t) = t_3 + t, & t \in [T_q - T_p, 0], \\
h(t) \; \mbox{increases}, & t \in [T_p, T_p +T'] \cup [T_q + T'',
T_q].
\end{cases}
$$
If $t \geq T_p + T'$, then inequality \sref{Pst1.5.1} implies that
$$
\dist(\phi(h(t), x_0), \phi(t - (T_p + T'), z_p)) < \ep/4
$$
and
$$
\dist(\phi(t - T_p, x_p), \phi(t - (T_p + T'), z_p)) < \ep/4.
$$
Hence, if $t \geq T_p + T'$, then
\begin{equation}\label{Pst1.7.1}
\dist(\phi(h(t), x_0), g(t)) < \ep/2.
\end{equation}
Inclusion \sref{Text23.1} implies that for $t \in [T_p, T_p + T']$
the inclusions $\phi(h(t), x_0), g(t) \in B(m, p)$ hold, and
inequality \sref{Pst1.7.1} holds for these $t$ as well.

A similar reasoning shows that inequality \sref{Pst1.7.1} holds for
$t \leq T_q$. If $t \in [T_q, T_p]$, then inequality \sref{Pst1.7.1}
follows from \sref{Pst1.6.5}. This completes the proof in case (P1).
\medskip

Case (P2): $x_p \in \Ws(p)$ and $x_q \notin \Wu(q)$. In this case,
Lemma 6 is replaced by the following statement.
\medskip

{\bf Lemma 7. } {\em There exists $m >0$ such that for any points
$$
y_p \in B(m, p) \cap \al, \quad y_q \in B(m, q) \cap \al, \quad z_q
\in W^s_{loc}(q, m)  \setminus \{q\},
$$
and a number $m_1 > 0$ there exists a trajectory of the vector field
$X$ that intersects successively the balls $B(m_1, z_q)$, $B(m_1,
y_q)$, and $B(m_1, y_p) \cap W^s_{loc}(p, m)$ as time grows.}
\medskip

The rest of the proof uses the same reasoning as in case (P1).
\medskip

Case (P3): $x_p \notin \Ws(p)$ and $x_q \in \Wu(q)$. This case is
similar to case (P2).
\medskip

Case (P4): $x_p \in \Ws(p)$ and $x_q \in \Wu(q)$. In this case, we
take $\al$ as the shadowing trajectory; the reparametrization is
constructed similarly to case (P1).
\medskip

Thus, to complete the consideration of case (S2), it remains to
prove Lemmas 6 and 7.

To prove Lemma 6, we first fix proper coordinates in small
neighborhoods of the points $p$ and $q$. Let us begin with the case
of the point $p$.

Taking a small neighborhood $V$ of the vector field $X^*$, we may
assume that the Jacobi matrix $J_p=DX(p)$ is as close to $J^*_p$ as
we want.

Thus, we assume that $p=0$ in coordinates $u_1 = (x_1, x_2)$, $u_2 =
(x_3, x_4)$, and $J_p=\diag(A_p,B_p)$, where
\begin{equation} \label{Text12.0.5}
A_p = \left(
      \begin{array}{cc}
        -\lam_1 & 0 \\
        0 & -\lam_2 \\
      \end{array}
    \right), \quad
B_p = \left(
        \begin{array}{cc}
          a_p & -b_p \\
          b_p & a_p \\
        \end{array}
      \right),
\end{equation}
and
\begin{equation}\label{Text4.5}
\lam_1, \lam_2, a_p,  b_p > 4g,
\end{equation}
where $g$ is a small positive number to be chosen later (and a
similar notation is used in $U_q$).

Then we can represent the field $X$ in a small neighborhood $U$ of
the point $p$ in the form
\begin{equation} \label{Text12.1}
 X(u_1, u_2) = \left(
                \begin{array}{cc}
                  A_p & 0 \\
                  0 & B_p \\
                \end{array}
              \right) \left(
                        \begin{array}{c}
                          u_1 \\
                          u_2 \\
                        \end{array}
              \right) +
                \left(
                  \begin{array}{c}
                    X_{12}(u_1, u_2) \\
                    X_{34}(u_1, u_2) \\
                  \end{array}
                \right),
\end{equation}
where
\begin{equation}\label{SM1.2.1}
 X_{12}, X_{34} \in \Cone, \; |X_{12}|_{\Cone}, |X_{34}|_{\Cone} < g,
\; X_{12}(0, 0) = X_{34}(0, 0) = (0, 0).
\end{equation}
Under these assumptions, $p=0$ is a hyperbolic rest point whose
two-dimensional unstable manifold in the neighborhood $U$ is given
by $u_2 = G(u_1)$, where $G:\RR^{2} \to \RR^2$, $G \in \Cone$. We
can find $g > 0$ such that if the functions $X_{12}$ and $X_{34}$
satisfy relations \sref{SM1.2.1}, then
\begin{equation}\label{SM1.3}
\|D G(u_1)\| < 1 \quad \mbox{while} \quad (u_1, G(u_1)) \in U.
\end{equation}
We introduce new coordinates in $U$ by $v(u_1, u_2) = (u_1, u_2 -
G(u_1))$ and use a smooth cut-off function to extend $v$ to a
$\Cone$ diffeomorphism $w$ of $M$ such that $w(x)=x$ outside a
larger neighborhood $U'$ of $p$. Denote by $Y$ the resulting vector
field in the new coordinates.
\medskip

{\bf Remark 3. } Note that $Y$ is continuous but not necessary
$\Cone$. Nevertheless, the following holds. Let $S_1$ and $S_2$ be
small smooth three-dimensional disks transverse to a trajectory of
$Y$ and let $f_Y$ be the corresponding Poincar\'e transformation
generated by the vector field $Y$. Consider smooth disks
$w^{-1}(S_1)$ and $w^{-1}(S_2)$ and let $f_X:w^{-1}(S_1) \to
w^{-1}(S_2)$ be the corresponding Poincar\'e transformation. Since
$f_X \in \Cone$ and $f_Y = w \circ f_X \circ w^{-1}$, we conclude
that $f_Y \in \Cone$. We will use this fact below.
\medskip

If $(v_1, v_2) = v(u_1, u_2)$, then
\begin{equation}\label{Text13.1}
u_1 = v_1, \quad u_2 = v_2 + G(v_1).
\end{equation}
Let $Y(v_1, v_2) = (Y_1(v_1, v_2), Y_2(v_1, v_2))$. Since the
surface $u_2 = G(u_1)$ is a local stable manifold of the rest point
$0$ of the field $X$, the surface $v_2 = 0$ is a local stable
manifold of the rest point $0$ of the vector field $Y$. Hence,
$$
Y_2(v_1, 0) = 0 \quad  \mbox{for} \quad (v_1, 0) \in v(U).
$$

{\bf Lemma 8. } {\em The inequalities
\begin{equation}\label{SM2.0}
|Y_2(v_1, v_2) - (Y_2(v_1, 0)+B_pv_2)| \leq 2 g |v_2|, \quad (v_1,
v_2) \in v(U),
\end{equation}
hold.}
\medskip

{\em Proof.} Substitute equalities \sref{Text13.1} into
\sref{Text12.1} to show that
\begin{multline*}
Y_2(v_1, v_2) = B_p(v_2 + G(v_1)) + X_{34}(v_1, v_2 + G(v_1)) - \\ -
D G(v_1) (A_p v_1 + X_{12}(v_1, v_2 + G(v_1))).
\end{multline*}
Relations \sref{SM1.2.1} and \sref{SM1.3} imply that
$$
|X_{34}(v_1, v_2 + G(v_1)) - X_{34}(v_1, G(v_1))| \leq g|v_2|
$$
and
$$
|D G(v_1) (A_p v_1 + X_{12}(v_1, v_2 + G(v_1))) - D G(v_1) (A_p v_1
+ X_{12}(v_1, G(v_1)))| \leq g|v_2|.
$$
Hence,
\begin{multline*}
|X_{34}(v_1, v_2 + G(v_1)) - X_{34}(v_1, G(v_1)) - \\ - ( D G(v_1)
(A_p v_1 + X_{12}(v_1, v_2 + G(v_1))) - D G(v_1) (A_p v_1 +
X_{12}(v_1, G(v_1))) )| \leq
\\ \leq 2 g |v_2|.
\end{multline*}
The left-hand side of the above inequality equals $|Y_2(v_1, v_2) -
(Y_2(v_1, 0)+B_p v_2)|$, which proves inequality \sref{SM2.0}.

Note that if $y_p, y_q, z_p, z_q$, and $m_1>0$ are fixed, then there
exists $m^* > 0$ such that if a trajectory $\beta^*$ of the vector
field $Y$ intersects successfully the balls $B(m^*, v(z_q))$,
$B(m^*, v(y_q))$, $B(m^*, v(y_p))$, and $B(m^*, v(z_p))$, then the
trajectory $w^{-1}(\beta^*)$ of $X$ has the property described in
Lemma 6.

Thus, it is enough to prove Lemma 6 for the vector field $Y$. Since
the mapping $w$ is smooth, the vector field $Y$ satisfies condition
\sref{Add3.2}.

To simplify presentation, denote $Y$ by $X$ and its flow by $\phi$.
In this notation, there exists a neighborhood $U_p$ of $p=0$ in
which
\begin{equation}\label{BL1.1}
X(x) =  \left(
                \begin{array}{cc}
                  A_p & 0 \\
                  0 & B_p \\
                \end{array}
              \right)  x + X_p(x),
\end{equation}
where $X_p \in \Cc^0$, and if $(x_1, x_2, x_3, x_4) \in U_p$, then
\begin{equation}\label{BL0.1}
|P^p_{34} X_p(x_1, x_2, x_3, x_4) | < 2g \max(|x_3|, |x_4|) \quad
\mbox{and}\quad P^p_{34} X_p(x_1, x_2, 0, 0) = 0
\end{equation}
(where we denote by $P^p_{34}$ the projection in $U_p$ to the plane
of variables $x_3, x_4$ parallel to the plane of variables $x_1,
x_2$). Conditions \sref{BL0.1} imply that the plane $x_3 = x_4 =0$
is a local stable manifold for the vector field $X$.

Introduce polar coordinates $r, \varphi$ in the plane of variables
$x_3, x_4$. In what follows (if otherwise is not stated explicitly),
we use coordinates $(x_1, x_2, r, \varphi)$. For $i \in \{1, 2, 3,
4, r, \varphi \}$, we denote by $P^p_i x$ the $i$th coordinate of a
point $x \in U_p$.

Since the surface $\Wu(p)$ is smooth and transverse to the plane
$x_3 = x_4 = 0$, there exist numbers $K > 0$ and $m_2
> 0$ such that if points $x \in W^u_{loc}(p, m_2)$ and
$y \in B(m_2, p)$ satisfy the equality $P^p_{34}x = P^p_{34} y$,
then
\begin{equation}\label{BLAdd1.3}
\dist(x, y) \leq K \dist(y, W^u_{loc}(p, m_2)).
\end{equation}
We reduce the neighborhood $U_p$ so that $U_p \subset B(m_2, p)$.
\medskip

{\bf Lemma 9. } {\em Let $x(t) = (x_1(t), x_2(t), r(t), \varphi(t))$
be a trajectory of the vector field $X$. The relations
\begin{equation}\label{BL3.1}
\ddt{r} \in ((a_p - 4g)r, (a_p+ 4g)r) \quad\mbox{and}\quad \ddt{\vp}
\in (b_p - 4g, b_p + 4g)
\end{equation}
hold while} $x(t) \in U_p$.
\medskip

{\em Proof. } Let $x_3(t) = P_3^p x(t)$ and $x_4(t) = P_4^p x(t)$.
Relations \sref{Text12.0.5}, \sref{BL1.1} and \sref{BL0.1} imply
that
\begin{equation}\label{BLA4.1.1}\notag
\ddt{x_3}(t) = a_p x_3(t) - b_p x_4(t) + \Delta_3(t)
\end{equation}
and
\begin{equation}\label{BLA4.1.2}\notag
\ddt{x_4}(t) = b_p x_3(t) + a_p x_4(t) + \Delta_4(t),
\end{equation}
where
\begin{equation}\label{BLA4.2}
|\Delta_3(t)|, |\Delta_4(t)| < 2 g r(t).
\end{equation}
Since $x_3(t) = r(t) \cos \varphi(t)$ and $x_4(t) = r(t)
\sin\varphi(t)$, we obtain the equalities
\begin{equation}\label{BLA4.3} \notag
r  \ddt{\varphi} = r b_p + \Delta_4(t) \cos \vp - \Delta_3(t) \sin
\vp
\end{equation}
and
$$
\ddt{r} = a_p r +\Delta_3(t) \cos \vp + \Delta_4(t) \sin \vp.
$$
Inequalities \sref{BLA4.2} imply that
\begin{equation}\label{BLA4.6}\notag
b_p - 4g < \ddt{\vp} < b_p + 4g
\end{equation}
and
\begin{equation}\label{BLA5.1}\notag
(a_p - 4g)r < \ddt{r} < (a_p + 4g)r,
\end{equation}
which proves our lemma.

A similar reasoning shows that there exists a neighborhood $U_q$ of
the point $q$ in which we can introduce (after a smooth change of
variables) coordinates $(y_1, y_2, y_3, y_4)$ (and the corresponding
polar coordinates $(r, \varphi)$ in the plane of variables $y_3,
y_4$) such that
$$
W^u_{loc}(q, m) \subset \{y_3 = y_4 = 0\}
$$
and for any trajectory $y(t) = (y_1(t), y_2(t), r(t), \varphi(t))$
of the vector field $X$, the relations
\begin{equation}\notag
\ddt{r} \in ((a_q - 4g)r, (a_q+ 4g)r) \quad\mbox{and}\quad
\ddt{\vp}\in (-b_q - 4g, -b_q + 4g)
\end{equation}
hold while $y(t) \in U_q$.
\medskip

Let us continue the proof of Lemma 6.

Let $S_p \subset U_p$ and $S_q \subset U_q$ be smooth
three-dimensional disks that are transverse to the vector field $X$
and contain the points  $y_p$ and $y_q$, respectively. Denote by $f:
S_q \to S_p$ the corresponding Poincar\'e transformation (generated
by the field $X$). We note that $f \in \Cone$ (see Remark 3) and
$f(y_q) = y_p$.

Consider the lines $l_p = S_p \cap W^s_{loc}(p, m)$ and $l_q = S_q
\cap W^u_{loc}(q, m)$ and unit vectors $e_p \in l_p$ and $e_q \in
l_q$. Let $P^p_{34}$ and $P^q_{34}$ be the projections to the planes
of variables $x_3,x_4$ and $y_3,y_4$ in the neighborhoods $U_p$ and
$U_q$, respectively. Relation \sref{Add3.2} implies that
\begin{equation} \label{Text14.1}
P^p_{34} D f(y_q)e_q \ne 0 \quad\mbox{and}\quad  P_{34}^q D
f^{-1}(y_p) e_p \ne 0.
\end{equation}

Take $m_3 \in (0, m_1)$ such that
$$
\phi(t, x) \in U_p, \quad x \in B(m_3, y_p),\; t \in (0, \tau_p(x)),
$$
and
$$
\phi(t, y) \in U_q, \quad y \in B(m_3, y_q),\; t \in (\tau_q(x), 0),
$$
where
$$
\tau_p(x) = \inf \{ t>0 : \; P_r^p(\phi(t, x)) \geq P_r^p z_p \},
$$
$$
\tau_q(x) = \sup \{ t< 0 : \; P_r^q(\phi(t, y)) \geq P_r^q z_q\},
$$
and $z_p,z_q$ are the points mentioned in Lemma 6.

Consider the surface $L_p \subset S_p$ defined by
$$
L_p = \{x+(y-y_p), \; x \in l_p, y \in f(l_q)\}.
$$
Let $L_q = f^{-1}L_p \subset S_q$. The surfaces $L_p$ and $L_q$ are
divided by the lines $l_p$ and $l_q$ into half-surfaces. Let $L_p^+$
and $L_q^+$ be any of these half-surfaces.

To any point $x \in L_p^+ \cap f(L_q^+)$ there correspond numbers
$r_p(x) = P^p_r x$ and ${r_q(x) = P^q_r f^{-1}(x)}$; consider the
mapping $w : L_p^+\cap f(L_q^+) \to \RR^2$ defined by $w(x) =
(r_p(x), r_q(x))$. We claim that there exists a neighborhood $U_L
\subset L_p^+ \cap f(L_q^+)$ of the point $y_p$ on which the mapping
$w$ is a homeomorphism onto its image.

Let $r_0$ and $\varphi_0$ be the polar coordinates of the vector
$P_{34}^p D f(y_q)e_q$. Relation \sref{Text14.1} implies that $r_0
\ne 0$. Hence, there exists a neighborhood $V_q$ of the point $y_q$
in $S_q$ such that if $y \in V_q$, then
\begin{equation}\label{BLAdd1.1}
P_r^p D f(y) e_q \in [r_0/2, 2r_0] \quad\mbox{and}\quad
P_{\varphi}^p D f(y) e_q \in [\varphi_0 -\pi/8, \varphi_0 + \pi/8].
\end{equation}

Take $c > 0$ such that $B(2c, y_q) \subset V_q$. Note that
$$
f(y_q + \delta e_q) = f(y_q) + \int_0^\delta D f(y_q + s e_q) e_q
\dd s, \quad \delta \in [0, c].
$$
Conditions \sref{BLAdd1.1} imply that
\begin{equation}\label{BLAdd1.2}
P^p_{\varphi} \left( f(y_q + \delta e_q) - f(y_q) \right) \in
[\varphi_0 -\frac{\pi}{8}, \varphi_0 + \frac{\pi}{8}], \; \delta \in
[0, c],
\end{equation}
and the mapping $Q_p(\delta): [0, c] \to \RR$ defined by
$Q_p(\delta) = P^p_{r} f(y_q + \delta e_q)$ is a homeomorphism onto
its image. Similarly (reducing $g$, if necessary), one can show that
if $x \in B(g, y_p)$, then the mapping  $Q_{q, x}(\delta): [0, g]
\to \RR$ defined by $Q_{q, x}(\delta) = P^q_{r} f^{-1}(x + \delta
e_p)$ is a homeomorphism onto its image.

Take $\delta_p, \delta_q \in [0, c]$ and let $x = \delta_p e_p +
f(y_q + \delta_q e_q)$. Then $r_p(x) = Q_p(\delta_q)$ and $r_q(x) =
Q_{q, f(y_q +\delta_q e_q)}(\delta_p)$. It follows that the mapping
$w$ is a homeomorphism onto its image. Indeed, if $g_1> 0$ is small
enough, then the mapping $w^{-1}(\xi, \eta) = (x(\xi), Q_{q,
x(\xi)}^{-1}(\eta))$, where $x(\xi) = f(y_q + Q_p^{-1}(\xi)e_q)$, is
uniquely defined and continuous for $(\xi, \eta) \in [0, g_1]\times
[0, g_1]$.

We reduce $m_3$ so that the following relations hold:
$$
m_3 < c, \quad B(m_3, y_p) \cap L^+_p \subset U_L,
\quad\mbox{and}\quad B(m_3, y_q) \cap L^+_q \subset f^{-1} U_L.
$$
Let us prove a statement which we use below.
\medskip

{\bf Lemma 10. } {\em For any $m_1> 0$ there exist numbers $r_1, r_2
\in (0, m_1)$ and $T_1, T_2 > 0$ with the following property: if
$\gam(s): [0, 1] \to L_p^+$ is a curve such that
\begin{equation}\label{BL4.0.5}
P^p_r \gam(0) = r_1, \quad P^p_r \gam(1) = r_2,
\end{equation}
and
\begin{equation}\label{BL4.1}
\gam(s) \in L_p^+ \cap B(m_2, y_p), \quad s \in [0, 1],
\end{equation}
then there exist numbers $\tau\in[T_2,T_1]$ and $s \in [0, 1]$ such
that}
$$
\phi(\tau,\gam(s))\in B(m_1,z_p).
$$

{\em Proof.} Let $r_p = P^p_r z_p$ and $\varphi_p = P^p_{\varphi}
z_p$. For $r> 0$, denote
$$
T_{\min}(r) = \frac{\log r_p - \log r}{a_p+4g} \quad\mbox{and}\quad
T_{\max}(r) = \frac{\log r_p - \log r}{a_p-4g}.
$$
Note that if $r < r_p$, then $T_{\max}(r) > T_{\min}(r)$ and that
$T_{\min}(r) \to \infty$ as $r \to 0$. Take $T > 0$ such that if
$\tau>T$, $x\in B(m_2, y_p)$, and
$$
\phi(t, x) \subset U_p,\quad t\in [0, \tau],
$$
then
\begin{equation}\label{BL4.2}
\dist(W^u_{loc}(p, m), \phi(\tau, x)) < \frac{m_1}{2K}.
\end{equation}
Take $r_1, r_2 \in (0, \min(m_2, r_p))$ such that
\begin{equation}\label{BL4.3}\notag
r_2> r_1, \quad T_{\min}(r_2) > T,
\end{equation}
and
\begin{equation}\label{BL5.1}
(b_p-4g)T_{\min}(r_1) - (b_p+4g)T_{\max}(r_2) > 4 \pi.
\end{equation}
Set $T_1 = T_{\max}(r_1)$ and $T_2 = T_{\min}(r_2)$. Since the
function $\gam(s)$ is continuous, inclusions \sref{BL3.1} and
inequalities \sref{Text4.5} imply that there exists a uniquely
defined continuous function ${\tau(s): [0, 1] \to \RR}$ such that
\begin{equation}\label{BL5.2}\notag
P^p_r \phi(\tau(s), \gam(s)) = r_p.
\end{equation}

It follows from inclusions \sref{BL3.1} and equalities
\sref{BL4.0.5} that
$$
\tau(0) \in [T_{\min}(r_1), T_{\max}(r_1)], \; \tau(1) \in
[T_{\min}(r_2), T_{\max}(r_2)], \; \tau(s) \in [T_2, T_1].
$$
Now we apply relations \sref{Text4.5}, \sref{BL3.1}, and
\sref{BLAdd1.2} to show that
$$
P^p_\varphi \phi(\tau(0), \gam(0)) \geq (b_p-4g)T_{\min}(r_1) +
\varphi_{0} - \pi/8
$$
and
$$
P^p_\varphi \phi(\tau(1), \gam(1)) \leq (b_p+4g)T_{\max}(r_2) +
\varphi_{0} + \pi/8.
$$
Since the function $\tau(s)$ is continuous, the above inequalities
and inequalities \sref{BL5.1} imply the existence of $s \in [0, 1]$
such that
$$
P^p_{\varphi} \phi(\tau(s), \gam(s))= \varphi_p \mod 2\pi.
$$
Hence, $P_{34}^p \phi(\tau(s),\gam(s)) = P_{34}^p z_p$. It follows
from this equality combined with relations \sref{BLAdd1.3},
\sref{BL4.2}, and the inequality $\tau(s) > T$  that $\phi(\tau(s),
\gam(s)) \in B(m_1/2, z_p)$, which proves Lemma 10.

Let $r_1, r_2 \in (0, m_2)$ and  $T_1, T_2 > 0$ be the numbers given
by Lemma 10. Consider the set
$$
A_p = \{\phi(t, x): \; t \in [-T_1, -T_2], x \in \Cl B(m_2/2, z_p)\}
\cap L_p^+.
$$
Note that $A_p$ is a closed set that intersects any curve $\gam(s)$
satisfying conditions \sref{BL4.0.5} and \sref{BL4.1}.

We apply a similar reasoning in the neighborhood $U_q$ to the vector
field $-X$ to show that there exist numbers ${r'_1, r'_2 \in (0,
m_2)}$ and ${T'_1, T'_2 > 0}$ such that the set
$$
A_q = \{\phi(t, x): \; t \in [T'_2, T'_1], x \in \Cl B(m_2/2, z_q)\}
\cap L_q^+
$$
is closed and intersects any curve $\gam(s): [0, 1] \to L_q^+ \cap
B(m_2, y_q)$ such that
\begin{equation}\label{Text19.1}\notag
P^q_r \gam(0) = r'_1 \quad\mbox{and}\quad P^q_r \gam(1) = r'_2.
\end{equation}

We claim that
\begin{equation}\label{BL6.1}
A_p \cap f(A_q) \ne \emptyset,
\end{equation}
which proves Lemma 6.

Consider the set $K \subset L_p^+ \cap f(L_q^+)$ bounded by the
curves $k_1 = {L_p^+ \cap \{P^p_r x = r_1\}}$, $ k_2 = L_p^+ \cap
\{P^p_r x = r_2\}$, $k'_1 = f(L_q^+ \cap \{P^q_r y = r'_1\})$, and
$k'_2 = {f(L_q^+ \cap \{P^q_r y = r'_2\})}$. Since $w(x)$ is a
homeomorphism, the set $K$ is homeomorphic to the square $[0, 1]
\times [0, 1]$.

The following statement was proved in \cite{PilSak}.
\medskip

{\bf Lemma 11. }{\em Introduce in the square $I =[0, 1]\times[0, 1]$
coordinates $(u, v)$. Assume that closed sets $A, B \subset I$ are
such that any curve inside $I$ that joins the segments $u = 0$ and
$u = 1$ intersects the set $A$ and any curve inside $I$ that joins
the segments
 $v = 0$ and $v = 1$ intersects the set $B$.
Then} $A \cap B \ne \emptyset$.

The set $A_p$ is closed. By Lemma 10, $A_p$ intersects any curve in
$K$ that joins the sides
 $k_1$ and $k_2$. Similarly, the set $A_q$ is closed and
intersects any curve that belongs to $f^{-1}(K)$  and joins the
sides $f^{-1}(k'_1)$ and $f^{-1}(k'_2)$. Thus, the set $f(A_q)$
intersects any curve in $K$ that joins the sides $k'_1$ and $k'_2$.
By Lemma 11 inequality \sref{BL6.1} holds. Lemma 6 is proved.
\medskip

{\em Proof of Lemma 7. } Similarly to the proof of Lemma 6, let us
consider the subspaces $L_p^+$ and  $L_q^+$ and a number  $m_2 \in
(0, m_1)$ and construct the set $A_q \subset L_q^+$. Note that the
set $f^{-1}(B(m_1, y_p) \cap \Ws(p) \cap L_p^+)$ contains a curve
that satisfies conditions (\ref{BL4.0.5}) and (\ref{BL4.1}). Hence,
$B(m_1, y_p) \cap \Ws(p) \cap f(A_q) \ne \emptyset$. For any point
in this intersection, its trajectory is the desired shadowing
trajectory.

\section{Appendix: Construction of the vector field $X^*$}

Consider two  2-dimensional spheres $M_1$ and $M_2$. Let us
introduce coordinates $(r_1, \varphi_1)$ and $(r_2, \varphi_2)$ on
$M_1$ and $M_2$, respectively, where $r_1, r_2 \in [-1, 1]$ and
$\varphi_1, \varphi_2 \in \RR/2\pi\Z$. We identify all points of the
form $(-1, \cdot)$ as well as points of the form $(1, \cdot)$.
Denote
$$
M_1^+ = \{(r_1, \varphi_1), \quad r_1 \geq 0\}\mbox{ and }  M_1^- =
\{(r_1, \varphi_1), \quad r_1 \leq 0\}.
$$
Consider a smooth vector field $X_1$ defined on $M_1^+$ such that
its trajectories $(r_1(t), \varphi_1(t))$ satisfy the following
conditions:
$$
\ddt r_1  = 1, \ddt \varphi_1 = 0, \quad r_1 = 0;
$$
$$
\ddt r_1 > 0, \quad r_1 >0;
$$
$$
\ddt r_1 = 0, \quad r_1 = 1.
$$
We also assume that, in proper local coordinates in a neighborhood
of the ``North Pole" $(1, \cdot)$ of the sphere $M_1$, the vector
field $X_1$ is linear, and
$$
\D X_1(1, \cdot) = \left(
                                 \begin{array}{cc}
                                   -2 & 0 \\
                                   0 & -1 \\
                                 \end{array}
                               \right).
$$
Thus, $(1, \cdot)$ is an attracting hyperbolic rest point of $X_1$,
and every trajectory of $X_1$ in $M_1^+$ tends to $(1,\cdot)$ as
time grows.

Consider a smooth vector field $X_2$ on $M_2$ such that its
nonwandering set $\Omega(X_2)$ consists of two rest points: a
hyperbolic attractor $s_2 = (0, \pi)$ and a hyperbolic repeller $u_2
= (0, 0)$. Assume that, in proper coordinates, the vector field
$X_2$ is linear in neighborhoods of $s_2$ and $u_2$, and
$$
\D X_2(s_2) = -\D X_2(u_2) = \left(\begin{array}{cc}
                                   -1 & 1 \\
                                   -1 & -1 \\
                                 \end{array}
                               \right).
$$

Consider the vector field $X^+$ defined on $M_1^+ \times M_2$ by the
following formula
$$
X^+(r_1, \varphi_1, r_2, \varphi_2) = (X_1(r_1, \varphi_1), r_1^2
X_2(r_2, \varphi_2)).
$$

Consider infinitely differentiable functions $g_1: M_1^+ \to \RR$,
$g_2, g_3: [-1, 1] \to [-1, 1]$, and $g_4: M_1^+ \to [0, 1]$
satisfying the following conditions:
$$
g_1(0, 0) = 0; \quad g_1(r_1, \varphi_1)\in (0, 2\pi), \quad (r_1,
\varphi_1) \ne 0,
$$
$$
g'_2(r_2) \in (0, 2), \quad r_2 \in [-1, 1];
$$
$$
g_2(0)< 0, \; g_2(-1) = -1, \; g_2(1) = 1;
$$
$$
g_3(r_2) = 2r_2 - g_2(r_2), \quad r_2 \in [-1, 1];
$$
$$
g_4(0, 0) = 1/2, \; \frac{\partial}{\partial \varphi_1} g_4(0,0) \ne
0.
$$
Note that the functions $g_2$ and $g_3$ are monotonically
increasing.

Consider a mapping $f^* : M_1^+ \times M_2 \to M_1^- \times M_2$
defined by the following formula:
$$
f^*(r_1, \varphi_1, r_2, \varphi_2) = (-r_1, \varphi_1, g_4(r_1,
\varphi_1)g_2(r_2) + (1 - g_4(r_1, \varphi_1))g_3(r_2), \varphi_2 +
g_1(r_1, \varphi_1)).
$$
Clearly, $f^*$ is surjective; the monotonicity of $g_2$ and $g_3$
implies that $f^*$ is a diffeomorphism.

Using the standard technique with a ``bump'' function, one can
construct a diffeomorphism $f : M_1^+ \times M_2 \to M_1^- \times
M_2$ such that, for small neighborhoods $U_1 \subset U_2$ of $(1,
\cdot, s_2)$, the following holds:
$$
f(x) = f^*(x), \quad x \notin U_2,
$$
and $f$ is linear in $U_1$.

Consider the set $l = \{r_1 = 0, r_2 = 0, \varphi_2 = 0\}$. Simple
calculations show that
\begin{equation}\label{S2.0}
f(l) \cap l = \{(0, 0, 0, 0)\},
\end{equation}
and the tangent vectors to $l$ and $f(l)$ at $(0, 0, 0, 0)$ are
parallel to the vectors $(0, 1, 0, 0)$ and $(0, 1, (g_2(0) -
g_3(0))\frac{\partial}{\partial \varphi_1}g_4(0, 0), \cdot)$,
respectively. Hence,
\begin{equation}
\label{S2.0.5} \dim(T_{(0, 0, 0, 0)}l \oplus T_{(0, 0, 0, 0)}f(l)) =
2.
\end{equation}

Define a vector field $X^-$ on $M_1^- \times M_2$ by the formula
$$
X^-(x) = -\D f(f^{-1}(x))X^+(f^{-1}(x))
$$
(and note that $x(t)$ is a trajectory of $X^+$ if and only if
$f(x(-t))$ is a trajectory of $X^-$).

Finally, we define the following vector field $X^*$ on $M_1 \times
M_2$:
$$
X^*(x) = \begin{cases} X^+(x), & \quad x \in M_1^+ \times M_2,\\
X^-(x), & \quad x \in M_1^- \times M_2
\end{cases}
$$

Let us check that the vector field $X^*$ is well-defined on the set
$\{r_1 =0\}$. Indeed, $X^+(0, \varphi_1, r_2, \varphi_2) = (1, 0, 0,
0)$ and $(\D f(0, \varphi_1, r_2, \varphi_2))^{-1} (1, 0, 0, 0) =
(-1, 0, 0, 0)$. It is easy to see that $\D X^+(0, \varphi_1, r_2,
\varphi_2) = \D X^-(0, \varphi_1, r_2, \varphi_2) = 0$. This implies
that $X \in \Cone$.

Let us prove that the vector field $X^*$ satisfies conditions (F1)
-- (F3). Let $(r_1(t), \varphi_1(t), r_2(t), \varphi_2(t))$ be a
trajectory of $X^*$. The following inequalities hold:
\begin{equation}
\label{S2.1} \ddt r_1 > 0, \quad r_1 \ne \pm 1.
\end{equation}
This implies the inclusion $\Omega(X^*) \subset \{r_1 = \pm 1\}$. By
the construction of $X^+$, $\Omega(X^*) \cap \{r_1 = 1\} = \{(1,
\cdot, s_2), (1, \cdot, u_2)\}$. Similarly, $\Omega(X^*) \cap \{r_1
= -1\} = \{f(1, \cdot, s_2), f(1, \cdot, u_2)\}$. Denote $s^* = (1,
\cdot, s_2)$, $p^* = (1, \cdot, u_2)$, $q^* = f(p)$, and $u^* =
f(s)$. Clearly, $s^*$, $u^*$, $p^*$, $q^*$ are hyperbolic rest
points, $s^*$ is an attractor, $u^*$ is a repeller, $\D X(p^*) =
J_p^*$, and $\D X(q^*) = J_q^*$. In addition, in small neighborhoods
of $p^*$ and $q^*$, the vector field $X^*$ is linear.

It is easy to see that
$$
\Ws(p^*) \cap \{r_1 = 1\} = \{p^*\}\mbox{ and } \Ws(p^*) \cap \{r_1
=-1\} = \emptyset.
$$
Inequality \sref{S2.1} implies that any trajectory in
$\Ws(p^*)\setminus\{p^*\}$ intersects the set $\{r_1 = 0\}$ at a
single point. The definition of $X^+$ implies that $\Ws(p^*) \cap
\{r_1 = 0\} = l$. Similarly, any trajectory in
$\Wu(q^*)\setminus\{q^*\}$ intersects $\{r_1 = 0\}$ at a single
point, and $\Wu(q^*) \cap \{r_1 = 0\} = f(l)$. It follows from
equality \sref{S2.0} that $\Ws(p^*) \cap \{r_1 = 0\} \cap \Wu(q^*)$
is a single point, and hence $\Ws(p^*) \cap \Wu(q^*)$ consists of a
single trajectory.

Inequality \sref{S2.1} implies condition \sref{Text16.5}, and
condition \sref{S2.0.5} implies \sref{dim3}.
\medskip

The authors are deeply grateful to the anonymous referee whose
remarks helped us to significantly improve the presentation.

\end{document}